\documentclass[10pt]{amsart}

\usepackage{fullpage}
\usepackage{amssymb}
\usepackage{amsmath}
\usepackage{amsxtra}
\usepackage{amscd}
\usepackage{graphicx}
\usepackage{amsfonts}
\usepackage{pb-diagram}
\usepackage{float}
\usepackage{enumerate}
\usepackage{dsfont}
\usepackage{multirow}
\usepackage{color}
\usepackage{rotating}
\usepackage{url}
\usepackage{enumitem}
\usepackage[bookmarks=false]{hyperref}
\definecolor{darkgreen}{rgb}{0,0.4,0.1}
\definecolor{darkpurple}{rgb}{0.7,0.0,0.4}
\hypersetup{
   colorlinks=true,       	
   linkcolor=darkpurple,     
   citecolor=darkgreen,    
   filecolor=magenta,      
   urlcolor=blue			
}

\numberwithin{equation}{section}
\newtheorem{thm}{Theorem}[section]

\newtheorem{prop}[thm]{Proposition}

\newtheorem{conj}[thm]{Conjecture}
\newtheorem{exmp}[thm]{Example}

\newtheorem{rem}[thm]{Remark}

\def\Q{{\mathbb Q}}

\def\Z{{\mathbb Z}}

\def\C{{\mathbb C}}

\def\Irr{{\mathrm{Irr}}}

\title{The BMM symmetrising trace conjecture for groups $G_4,\,G_5,\,G_6,\,G_7,\,G_8$}
\author{Christina Boura}
\address{Laboratoire de Math\'ematiques UVSQ, B\^atiment Descartes, 45 avenue des \'Etats-Unis,  78035 Versailles cedex, France.}
\email{christina.boura@uvsq.fr}

\author{Eirini Chavli}
\address{Institut f\"ur Algebra und Zahlentheorie, Universit\"at Stuttgart, Pfaffenwaldring 57, 
70569 Stuttgart, Germany.}
\email{eirini.chavli@mathematik.uni-stuttgart.de}

\author{Maria Chlouveraki}
\address{Laboratoire de Math\'ematiques UVSQ, B\^atiment Fermat, 45 avenue des \'Etats-Unis,  78035 Versailles cedex, France.}
\email{maria.chlouveraki@uvsq.fr}

\author{Konstantinos Karvounis}
\address{Institut f\"ur Mathematik, 
Universit\"at Z\"urich, Winterthurerstrasse 190, CH-8057 Z\"urich, Switzerland.}\email{konstantinos.karvounis@math.uzh.ch}


\begin{document}

\maketitle

\begin{abstract}
We prove the BMM symmetrising trace conjecture for the exceptional  irreducible complex reflection groups $G_4,\,G_5,\,G_6,\,G_7,\,G_8$ using a combination of algorithms programmed in different  languages (\texttt{C++}, SAGE, GAP3, \emph{Mathematica}). Our proof depends on the choice of a  suitable  basis for the generic Hecke algebra associated with each group.
\end{abstract}

\section{Introduction}\label{Introduction}

Exactly twenty years ago, Brou\'e, Malle and Rouquier published their seminal paper \cite{BMR}  in  which they 
associated to every complex
reflection  group  two  objects  which  were  classically  associated  to
real reflection groups: a braid group and a Hecke algebra. 
Their work was further motivated by the theory, developed together with Michel  \cite{BMM}, 
that complex reflection groups could play the role of Weyl groups of
objects that generalise finite reductive groups, named ``Spetses''.
The four of them advocated that several nice properties of braid groups and Hecke algebras generalise from the real to the complex case, culminating in two main conjectures as far as the Hecke algebras are concerned: the ``freeness conjecture'', stated in \cite{BMR}, and the ``symmetrising trace conjecture'', stated in \cite{BMM}. 
The former states that the Hecke algebra associated with a complex reflection group is a free module of rank equal to the order of the group. The latter states that the Hecke algebra possesses a ``canonical'' symmetrising trace.
The first conjecture was a known fact  for all real reflection groups and for the groups of the infinite series $G(l, p, n)$. The second conjecture was also known to hold for all real reflection groups, and was  partially solved for the groups $G(l,p,n)$: for these groups, a symmetrising trace does exist, but it is still unclear whether it satisfies all canonicality conditions.  However, even though these two conjectures are the cornerstones in the study  of subjects that have flourished in the past twenty years, namely
 the representation theory of Hecke algebras associated with complex reflection groups, as well as of other related structures such as Cherednik algebras,  
 there had been until a few years ago little progress in their proof for the exceptional complex reflection groups (with the exception of $G_4$, which is the smallest exceptional complex reflection group).

In the past five years there have been outstanding developments regarding the BMR freeness conjecture. Taking into account some preprints that appeared in 2017, the BMR freeness conjecture is now proved for all the exceptional complex reflection groups --- an exhaustive list of references can be found in $\text{Section \ref{Conjectures}}$.
However, the BMM symmetrising trace conjecture has remained widely open, except for the groups $G_4$, $G_{12}$, $G_{22}$ and $G_{24}$, for which the conjecture was established by Malle and Michel in
\cite{MM10} (the case of $G_4$
was later independently checked by Marin and Wagner in \cite{Mar46}). The two main difficulties for tackling this conjecture are the following:
\begin{itemize}
\item[(1)] Up to recently, we did not possess bases for the generic Hecke algebras of exceptional complex reflection groups (the BMR freeness conjecture was still open). \smallbreak
\item[(2)] The standard programming language used when working with complex reflection groups and their Hecke algebras, GAP3, is not efficient  when we need to make complicated calculations with multivariable polynomials.  
\end{itemize}

In order to understand why these two issues made the verification of the BMM symmetrising trace conjecture so difficult until now, we have to explain what  a canonical symmetrising trace is. Let $W$ be a complex reflection group and let $\mathcal{H}(W)$ be the generic Hecke algebra associated with $W$ defined over a Laurent polynomial ring $R$ in several indeterminates. Let $\mathcal{B}$ be a basis for $\mathcal{H}(W)$ as an $R$-module. A symmetrising trace is a linear map $\tau : \mathcal{H}(W) \rightarrow R$ such that the matrix $A:=(\tau(bb'))_{b,b' \in \mathcal{B}}$ is symmetric and its determinant is a unit in $R$. The symmetrising trace $\tau$ is ``canonical'' when it satisfies two further conditions: one of them is that it specialises to the canonical symmetrising trace of the group algebra of $W$ (that is, $\tau(w)=\delta_{1w}$ for all $w \in W$) and the other one concerns the behaviour of its values when the indeterminates of the Laurent polynomial ring are inverted. These conditions imply that the trace $\tau$ is unique  \cite[2.1]{BMM}. Malle and Michel \cite{MM10} even conjectured that there exists a basis $\mathcal{B}$ for $\mathcal{H}(W)$, each of its elements lifting an element of the group $W$, such that $1 \in \mathcal{B}$
and  $\tau(b)=\delta_{1b}$ for all $b \in \mathcal{B}$. Then $\tau$ is the map that sends any element $h \in \mathcal{H}(W)$ to the coefficient of $1$ when $h$ is expressed as an $R$-linear combination of the elements of $\mathcal{B}$.
Unfortunately, this property cannot hold for every choice of basis for $\mathcal{H}(W)$. A very interesting fact is that, assuming that the canonical symmetrising trace exists, it can be expressed as a linear combination with non-zero coefficients of the irreducible characters of the algebra $\mathcal{H}(W)$  over its splitting field. The inverses of the coefficients are the Laurent polynomials known as ``Schur elements''. The Schur elements control much of the modular representation theory of the Hecke algebra and they have already been determined by Malle  \cite{Ma2, Ma5} for all non-real exceptional complex reflection groups at around the same time that the BMM symmetrising trace conjecture was formulated. 

It is now clear that a lack of basis for the Hecke algebra makes the verification of the BMM symmetrising trace conjecture impossible. However, with the recent progress on the BMR freeness conjecture, we now have bases for all exceptional complex reflection groups. In particular, the second author of this paper has constructed explicit bases for the generic Hecke algebras associated with the groups $G_4,\ldots,G_{16}$ in \cite{Ch17, Ch18}.

Nevertheless, even if we have  a basis, and even if we have a good basis so that $\tau$ can be defined as $1$ on $1$ and $0$ elsewhere, there is a considerable amount of calculations to be made. We need to have a computer program that can calculate all the entries of the $|W| \times |W|$ matrix $A$ and then compute its determinant. The GAP3 package CHEVIE \cite{chevie1,chevie2} contains a lot of important data on the representation theory of Hecke algebras (such as the irreducible characters and the Schur elements), but it is inefficient when dealing with multivariable Laurent polynomials. Moreover, it contains no method to express an element $h \in \mathcal{H}(W)$ as a linear combination of the elements of a given basis.

In this article, we overcome all the difficulties above, and we prove: 

\begin{thm}
The BMM symmetrising trace conjecture holds for the exceptional irreducible complex reflection groups $G_4$, $G_5$, $G_6$, $G_7$ and $G_8$. 
\end{thm}

The orders of these groups are respectively $24$, $72$, $48$, $144$ and $96$.
Our first step was to find a suitable basis $(\mathcal{B}_n)_{n=4,\ldots,8}$  for each generic Hecke algebra   $(\mathcal{H}(G_n))_{n=4,\ldots,8}$ (with $1\in \mathcal{B}_n)$ so that we could define the linear map $\tau$ on $\mathcal{H}(G_n)$ by setting $\tau(b):=\delta_{1b}$ for all $b \in \mathcal{B}_n$. It turns out that the bases given by the second author in \cite{Ch17} were perfect for $G_5$, $G_6$ and $G_8$, but we had to change the ones for $G_4$ and $G_7$. Especially for $G_7$, the change was not trivial as the reader will see in \S \ref{G7}. 
We  then created a program in the language \texttt{C++} which would write any product $(bb')_{b,b' \in \mathcal{B}_n}$ as a linear combination of elements of the basis, with the coefficient of $1$ being the corresponding entry of the matrix $A$. With the exception of $G_4$, this program was very time-consuming, so we came up with an elaborate algorithm which makes use of the inductive nature of the basis $\mathcal{B}_n$ in order to fill in the entries of the matrix $A$ row-by-row. Using inputs from the \texttt{C++} program, we created a second program in SAGE \cite{sagemath} which produced very quickly the matrix $A$ and its determinant. We verified that the matrix $A$ is symmetric and that its determinant is a unit in the ring of definition of  $\mathcal{H}(G_n)$ for all $n=4,\ldots,8$. Thus, we obtained that $\tau$ is a symmetrising trace on $\mathcal{H}(G_n)$ for all $n=4,\ldots,8$. We also checked the canonicality conditions, thus concluding that the BMM symmetrising trace conjecture holds for $G_4$, $G_5$, $G_6$, $G_7$ and $G_8$.

The paper is organised as follows: Sections  \ref{sym tr} and \ref{Conjectures} are preliminaries; the former contains information on symmetrising traces in general, while the latter introduces Hecke algebras and the conjectures about them. Section \ref{Algorithm} is the core of our paper. We present the outline of the algorithm for the \texttt{C++} program and the detailed algorithm for the SAGE program (Subsections \ref{C++} and \ref{SAGE} respectively). This is where we also give the bases $\mathcal{B}_n$ for all $n=4,\ldots,8$. In Subsection \ref{GAP3}, we briefly discuss an attempt of using GAP3 for our purposes, which worked only for the two smallest groups $G_4$ and $G_6$, and even then we had to transfer and run our program on \emph{Mathematica} because it is faster and deals better with multivariable polynomials. Finally, in Subsection \ref{extra con section}, we show how we verify the canonicality conditions for the symmetrising trace. In $\text{Section \ref{Results}}$, we summarise our results on the BMM symmetrising trace conjecture for groups $G_4$, $G_5$, $G_6$, $G_7$, $G_8$ and we explain how we also obtain a verification of the BMR freeness conjecture in the process. At the end of the paper, there is an appendix which contains some examples and calculations that we decided not to include in the main body of the paper for coherence reasons.

$ $\\
{\bf Acknowledgements:} We  are grateful to Jean Michel and Ivan Marin  for  enlightening discussions about the BMM symmetrising trace conjecture. We would also like to thank Gunter Malle and Meinolf Geck for their helpful comments and suggestions.

\section{Symmetrising traces}\label{sym tr}

Let $R$ be a commutative integral domain 
and let $A$ be an $R$-algebra, free and finitely generated as an $R$-module. If $R'$ is a commutative integral domain containing $R$, we will write $R'A$ for $R' \otimes_R A$. We will denote by $\mathrm{Irr}(R'A)$ the set of irreducible representations of $R'A$ and by $Z(R'A)$ the centre of $R'A$.

A {\em symmetrising trace} on the algebra $A$ is a linear map $\tau : A \rightarrow R$ such that:
\begin{itemize}
\item[(1)] $\tau(ab)=\tau(ba)$ for all $a,b \in A$,  and \smallbreak
\item[(2)]the bilinear form $A \times A \rightarrow R,\,\, (a,b) \mapsto \tau(ab)$ is non-degenerate.
\end{itemize}
If there exists a symmetrising trace on $A$, we say that $A$ is a {\em symmetric} algebra.

Let $\mathcal{B}$ be an $R$-basis of $A$. Then $\tau$ is a symmetrising trace on $A$ if and only if the matrix $(\tau(bb'))_{b,b' \in \mathcal{B}}$ is symmetric and its determinant is a unit in $R$.

\begin{exmp}\label{group}
{\rm Let $G$ be a finite group. The linear map $\tau: \Z[G] \rightarrow \Z$ defined by $\tau(1)=1$ and $\tau(g)=0$ for all $g \in G \setminus \{1\}$ is a symmetrising trace on $\Z[G]$; it is called the {\em canonical symmetrising trace} on $\Z[G]$.}
\end{exmp}

Suppose that there exists a symmetrising trace $\tau$ on $A$ and
let $K$ be a field containing $R$, such that the algebra $KA$ is split.
The map $\tau$ can be extended to $KA$ by extension of scalars. 
Let $E \in \Irr(KA)$ with  character $\chi_E$. 
We  have \cite[Lemma 7.1.7]{gepf}:
$$
\chi_E^\vee := \sum_{b \in \mathcal{B}} \chi_E(b)\, b^\vee \in Z(KA) \ ,
$$
where $(b^\vee)_{b\in\mathcal{B}}$ denotes the dual basis of $A$ with  respect to $\tau$,  that is, $\tau(b^\vee b') = \delta_{bb'}$. 
Due to  Schur's lemma,  any $z \in Z(KA)$ acts as a scalar on $E$; we denote this scalar by $\omega_E(z)$.
The $K$-algebra homomorphism $\omega_E:Z(KA) \rightarrow K$ is the {\em central character} associated with $E$. We define
$$s_E:=  \omega_{E}(\chi_E^\vee )$$
to be the {\em Schur element} associated with $E$.  We have $s_E \in R_K$, where $R_K$ denotes the integral closure of $R$ in $K$  \cite[Proposition 7.3.9]{gepf}. 
Moreover, note that the element $s_E$ satisfies:
$$
s_E \,\chi_E(1) =\sum_{b \in \mathcal{B}} \chi_E(b)\, \chi_E(b^\vee) \ .
$$

\begin{exmp}{\rm Let $G$ be a finite group and let $\tau$ be the canonical symmetrising trace on $A:=\Z[G]$. 
The set $\{g\}_{g \in G}$  forms a basis of $A$ over $\Z$, with $\{g^{-1}\}_{g \in G}$ the dual basis of $A$ with respect to $\tau$.
If $K$ is an algebraically closed field of characteristic $0$, then $KA$ is a split semisimple algebra and $s_E = |G|/\chi_E(1) \in \Q$ for all $E \in \Irr(KA)$. Because of the integrality of the Schur elements, we must have $|G|/\chi_E(1) \in  \Z_K \cap \Q = \Z$ for all $E \in \Irr(KA)$. Thus, we have also shown that $\chi_E(1)$ divides $|G|$.}
\end{exmp}

Now, the algebra $KA$ is semisimple if and only if $s_E \neq 0$ for all $E \in \Irr(KA)$. 
If this is the case, we  have:
\begin{equation}\label{sym form}
\tau = \sum_{E \in \Irr(KA)}\frac{1}{s_E}\chi_E.
\end{equation}
The above result is due to Curtis and Reiner \cite{CuRe}, but we follow the exposition in  
 \cite[Theorem 7.2.6]{gepf}.

\section{Conjectures}\label{Conjectures}

Let $V$ be a finite dimensional complex vector space. A \emph{pseudo-reflection} is a non-trivial element $s \in  \mathrm{GL}(V)$ that fixes a hyperplane pointwise,
that is, ${\rm dim}({\rm Ker}(s - {\rm id}_V))={\rm dim}(V)-1$. The hyperplane ${\rm Ker}(s - {\rm id}_V)$ is the \emph{reflecting hyperplane} of $s$.
A {\em complex reflection group} is a finite subgroup of $\mathrm{GL}(V)$ generated by pseudo-reflections. The classification of (irreducible) complex reflection groups is due to Shephard and Todd \cite{ShTo}:

\begin{thm}\label{ShToClas} Let $W \subset \mathrm{GL}(V)$ be an irreducible complex
reflection group (i.e., $W$ acts irreducibly on $V$). Then one of
the following assertions is true:
\begin{itemize}
  \item There exists a positive integer $n$ such that $(W,V) \cong (\mathfrak{S}_n, \C^{n-1})$.\smallbreak
  \item There exist positive integers $l,p,n$ with $l/p \in \Z$ and $l > 1$ such
  that $(W,V) \cong (G(l,p,n),\C^n)$, where $G(l,p,n)$ is the group of all 
  $n \times n$ monomial matrices whose non-zero entries are ${l}$-th roots of unity, while the product of all non-zero
  entries is an $(l/p)$-th root of unity. \smallbreak
  \item $(W,V)$ is isomorphic to one of the 34 exceptional groups
  $G_n$ $(n=4,\ldots,37)$.
\end{itemize}
\end{thm}

If $W \subset \mathrm{GL}(V)$ is an irreducible complex
reflection group, then the dimension of $V$ is called the {\em rank} of $W$.
We have ${\rm rank}(\mathfrak{S}_n)=n-1$,  ${\rm rank}(G(l,p,n))=n$ for $l>1$ and
${\rm rank}(G_n) \in \{2,3,\ldots,8\}$ for $n=4,\ldots,37$.
In particular, we have ${\rm rank}(G_n)=2$ for $n=4,\ldots,22$.
Furthermore, Benard \cite{Ben} and Bessis \cite{Bes1} have proved (using a case-by-case
analysis) that the field $K$ generated by the traces on $V$ of all the elements of $W$ is a splitting field for $W$. The field $K$  is called the \emph{field of
definition} of $W$.  If $K \subseteq \mathbb{R}$, then $W$ is a finite Coxeter group, and
 if $K=\mathbb{Q}$, then $W$ is a Weyl group.

\begin{rem}
	{\rm We have $G(1,1,n) \cong \mathfrak{S}_n$, $G(2,1,n) \cong B_n$,  $G(2,2,n) \cong D_n$, $G(m,m,2) \cong I_2(m)$, 
		$G_{23} \cong H_3$,  $G_{28}  \cong  F_4$, $G_{30}  \cong H_4$, $G_{35}  \cong  E_6$, $G_{36}  \cong  E_7$, $G_{37}  \cong E_8$.       
	}
\end{rem}

From now on, let $W$ be an irreducible complex reflection group. Let $\mathcal{A}$ be the set of reflecting hyperplanes of $W$ and let $V^{\textrm{reg}}:= V\setminus \bigcup_{H\in \mathcal{A}} H$. 
We set $P(W):= \pi_1(V^{\textrm{reg}}, x_0)$
and $B(W) := \pi_1(V^{\textrm{reg}}/W, x_0)$, where $x_0 \in V^{\textrm{reg}}$ is some fixed basepoint.
The group $P(W)$ is the {\em pure braid group} of $W$ and the group $B(W)$ is the {\em braid group} of $W$. It is known by \cite[Theorem 12.8]{Bes2} that
the centre of $B(W)$ is cyclic, generated by some element $\boldsymbol{\rm z}$. We set $\boldsymbol{\pi}: = \boldsymbol{\rm z}^{|Z(W)|} \in P(W)$, where $Z(W)$ denotes the centre of $W$.

For every orbit $\mathcal{C}$ of the action of $W$ on $\mathcal{A}$, let
$e_{\mathcal{C}}$ be the common order of the subgroups $W_H$, where $H$
is any element of $\mathcal{C}$ and $W_H$ is the pointwise stabiliser of $H$. Note that $W_H$ is cyclic, for all $H \in \mathcal{A}$. Let 
$\mathbb{Z}[\textbf{u},\textbf{u}^{-1}]$ denote the Laurent polynomial ring
in a set of indeterminates $\textbf{u}=(u_{\mathcal{C},j})_{(\mathcal{C} \in
\mathcal{A}/W)(0\leq j \leq e_{\mathcal{C}}-1)}$. The \emph{generic
Hecke algebra} \index{generic Hecke algebra} $\mathcal{H}(W)$ of $W$ is  the quotient of the group
algebra $\mathbb{Z}[\textbf{u},\textbf{u}^{-1}][B_W]$ by the ideal
generated by the elements of the form
$$(s-u_{\mathcal{C},0})(s-u_{\mathcal{C},1}) \cdots (s-u_{\mathcal{C},e_{\mathcal{C}}-1}),$$
where $\mathcal{C}$ runs over the set $\mathcal{A}/W$ and
$s$ runs over the set of monodromy generators around 
the images in $V^{\textrm{reg}}/W$ of 
the elements of $\mathcal{C}$ (see \cite[\S 2]{BMR} for their definition). 

\begin{exmp}\label{Example of G4} \rm
The generic Hecke algebra of $G_4$ is the $\Z[u_0^{\pm 1},u_1^{\pm 1},u_2^{\pm 1}]$-algebra (note that $s$ and $t$ are conjugate):
$$\mathcal{H}(G_4)=\left\langle s,\,t \,\,\,\left|\,\,\, sts=tst,\,\,\, \prod_{i=0}^2(s-u_i) =  \prod_{i=0}^2(t-u_i) =0 \right\rangle\right..$$
We can rewrite the presentation of $\mathcal{H}(G_4)$ as follows (with suitable elements $a,b,c$):
$$\mathcal{H}(G_4)=\left\langle s,\,t \,\,\,\left|\,\,\, sts=tst,\,\,\, s^3 = as^2+bs+c,\,\,\, t^3 = at^2+bt+c\right\rangle\right..$$
Note that, in the above presentation, only $c$ is a unit in the ring of definition of $\mathcal{H}(G_4)$.
Moreover, in the presentation of $\mathcal{H}(G_4)$ we have two kind of relations: the \emph{braid relation} $sts=tst$, and the \emph{positive Hecke relations}  $s^3 = as^2+bs+c$ and $t^3 = at^2+bt+c$. If we multiply the last two with $s^{-1}$ and $t^{-1}$ respectively, we obtain the \emph{inverse Hecke relations}  $s^{-1} =c^{-1}s^2-ac^{-1}s-bc^{-1}$ and $t^{-1} =c^{-1}t^2-ac^{-1}t-bc^{-1}$.
\end{exmp}

This definition of generic Hecke algebras of complex reflection groups is due to Brou\'e, Malle and Rouquier \cite{BMR},  who, together with Michel, also stated two fundamental conjectures about their structure: the ``freeness conjecture'' 
\cite[\S 4]{BMR}
and the ``symmetrising trace conjecture'' \cite[2.1]{BMM} -- both known facts in the case of real reflection groups \cite[IV, \S 2]{Bou05}.

\begin{conj}\label{BMR free}\
\textbf{``The BMR freeness conjecture''} The algebra $\mathcal{H}(W)$ is a free {\rm $\mathbb{Z}[\textbf{u},\textbf{u}^{-1}]$}-module of rank  $|W|$.
\end{conj}

 The following result   reduces the search for a basis to a search for a spanning set with the correct number of elements. For its proof, the reader may refer to \cite[Proof of Theorem 4.24]{BMR} or \cite[Proposition 2.4]{Mar43}.
 
 \begin{thm}\label{mention}
 If $\mathcal{H}(W)$ is generated as a {\rm $\mathbb{Z}[\textbf{u},\textbf{u}^{-1}]$}-module by  $|W|$ elements, then the BMR freeness conjecture holds.
 \end{thm}

Taking into account some very recent results, the BMR freeness conjecture is now a theorem. More precisely, it has been proved for: 
\begin{itemize}
\item the complex reflection groups of the infinite series $G(l,p,n)$ by \cite{ArKo, BM, Ar};\smallbreak
\item the group $G_4$ by \cite{BM, Fun, Mar41}; \smallbreak
\item the group $G_{12}$ by  \cite{MaPf}; \smallbreak
\item the groups $G_4,\ldots, G_{16}$ by \cite{Ch17, Ch18}; \smallbreak
\item the groups $G_{17},\,G_{18},\,G_{19}$ by \cite{Tsu} (with a computer method applicable to all rank $2$ groups);\smallbreak
\item the groups $G_{20}, \,G_{21}$ by \cite{MarNew};\smallbreak
\item the groups $G_{22},\ldots, G_{37}$ by \cite{Mar41, Mar43, MaPf}. \smallbreak
\end{itemize}

Malle \cite[5.2]{Ma4} has shown that, given that the BMR freeness conjecture holds, we can always find $N_W \in \Z_{>0}$ such that,  if we take
\begin{equation}\label{NW}
v_{\mathcal{C},j}^{N_W}:= \zeta_{e_\mathcal{C}}^{-j}u_{\mathcal{C},j} ,
\end{equation}
where $\textbf{{v}}:=(v_{\mathcal{C},j})_{(\mathcal{C} \in
\mathcal{A}/W)(0\leq j \leq e_{\mathcal{C}}-1)}$ is a set of indeterminates and
$\zeta_{e_\mathcal{C}}:=\exp(2\pi i/e_\mathcal{C})$, then the 
$K(\textbf{{v}})$-algebra
$K(\textbf{{v}})\mathcal{H}(W)$ is split semisimple.
 Taking $N_W$ to be the number of roots of unity in $K$ works every time, but sometimes it is enough to take $N_W$ to be even as small as $1$ (for example, if $W=G(l,1,n)$ or $W=G_4$).
Following Tits's deformation theorem,  the specialisation $v_{\mathcal{C},j}\mapsto 1$ induces a
bijection between
$\mathrm{Irr}(K(\textbf{v})\mathcal{H}(W))$ and  $\mathrm{Irr}(W)$. 

\begin{conj}\label{BMM sym}\
\textbf{``The BMM symmetrising trace conjecture''} There exists a (canonical) symmetrising trace $\tau$ on $\mathcal{H}(W)$ that satisfies the following two conditions:
\begin{enumerate}
\item $\tau$ specialises to the canonical symmetrising trace on the group algebra of $W$ when  $u_{\mathcal{C},j}\mapsto \zeta_{e_\mathcal{C}}^{j}$. \smallbreak 
\item  $\tau$ satisfies
\begin{equation}\label{extra con}
\tau(T_{\beta^{-1}})^* =\frac{\tau(T_{\beta\boldsymbol{\pi}})}{\tau(T_{\boldsymbol{\pi}})}, \quad \text{ for all } \beta \in B(W),
\end{equation}
where $\beta \mapsto T_{\beta}$ denotes the natural surjection $B(W) \rightarrow \mathcal{H}(W)$ and
$x \mapsto x^*$ the automorphism of  $\mathbb{Z}[\textbf{\em u},\textbf{\em u}^{-1}]$ given by $\textbf{\em u} \mapsto \textbf{\em u}^{-1}$.
\smallbreak
\end{enumerate}
\end{conj}

A symmetrising trace that satisfies Condition (1) is known to exist for the complex reflection groups of the infinite series $G(l,p,n)$ \cite{BreMa, MM98}; it is still unclear though whether this trace satisfies Condition (2).
On the other hand, the only (non-Coxeter) exceptional groups for which the trace conjecture has been proved so far
are $G_4$, $G_{12}$, $G_{22}$ and $G_{24}$
\cite{MM10} (the case of $G_4$
was later independently checked in \cite{Mar46}). In this article, we will prove it for the groups $G_4,\ldots,G_8$.

Note that if the BMR freeness conjecture holds, then $\tau$ is unique \cite[2.1]{BMM}. Moreover, following \eqref{sym form}, $\tau$ can be expressed as a linear combination of the irreducible characters 
$(\chi_E)_{E \in \Irr(W)}$ of the split semisimple algebra $K(\textbf{{v}})\mathcal{H}(W)$ as follows:
\begin{equation}\label{needed}
\tau = \sum_{E \in \Irr(W)}\frac{1}{s_E}\chi_E,
\end{equation}
where $s_E \in \Z_K[\textbf{{v}},\textbf{{v}}^{-1}]$ denotes the Schur element of $K(\textbf{{v}})\mathcal{H}(W)$  associated with $E \in \Irr(W)$. The Schur elements have been explicitly calculated for all complex reflection groups. In particular, for the non-real exceptional complex reflection groups, they have been completely determined by Malle \cite{Ma2, Ma5}.

All of our computational verifications in this article will depend on  the choice of a  suitable  basis for the generic Hecke algebra. Malle and Michel have stated the following conjecture  \cite[Conjecture 2.6]{MM10}, which we will also prove  for groups $G_4,\ldots,G_8$:

\begin{conj}\textbf{``The lifting conjecture''}
There exists a section (that is, a map admitting a left inverse) $W \rightarrow \boldsymbol{W}\subset B(W)$, $w \mapsto \boldsymbol{w}$ of $W$ in $B(W)$ 
such that $1 \in \boldsymbol{W}$, and such that for any $\boldsymbol{w} \in \boldsymbol{W}$ we have $\tau(T_{\boldsymbol{w}}) = \delta_{1\boldsymbol{w}}$.
\end{conj}

If the above conjecture holds, then Condition (1) of Conjecture \ref{BMM sym} is obviously satisfied. If further the elements $\{T_{\boldsymbol{w}}\,|\,\boldsymbol{w} \in \boldsymbol{W}\}$ form a
$\Z[\textbf{{u}},\textbf{{u}}^{-1}]$-basis of $\mathcal{H}(W)$, then, by \cite[Proposition 2.7]{MM10}, Condition (2) of Conjecture \ref{BMM sym} is equivalent to:
\begin{equation}\label{extra con2}
\tau(T_{\boldsymbol{w}^{-1}\boldsymbol{\pi}})=0, \quad \text{ for all } \boldsymbol{w} \in \boldsymbol{W} \setminus \{1\}.
\end{equation}
Using \eqref{needed}, this is in turn equivalent to:
\begin{equation}\label{extra con3}
\sum_{E \in \Irr(W)}\frac{\omega_E(T_{\boldsymbol{\pi}})}{s_E}\chi_E(T_{\boldsymbol{w}^{-1}})=0, \quad \text{ for all } \boldsymbol{w} \in \boldsymbol{W} \setminus \{1\}.
\end{equation}
An explicit formula for $\omega_E(T_{\boldsymbol{\pi}})$ is given by \cite[Equation (1.22)]{BMM}.

\section{Algorithms}\label{Algorithm}
In \cite{Ch17}, the second author proved the BMR freeness conjecture for the exceptional groups $G_4,\ldots,G_{15}$ by providing an explicit basis for the generic Hecke algebra in each case. 
In this section, we give a basis $(\mathcal{B}_n)_{n=4,\ldots,8}$  for each generic Hecke algebra   $(\mathcal{H}(G_n))_{n=4,\ldots,8}$ and we define a linear map $\tau$ on $\mathcal{H}(G_n)$ by setting $\tau(b):=\delta_{1b}$ for all $b \in \mathcal{B}_n$.  Our aim is to calculate the matrix
$A:=(\tau(bb'))_{b,b'\in \mathcal{B}_n}$, and to show that it is symmetric and has determinant invertible in
$\Z[\textbf{{u}},\textbf{{u}}^{-1}]$. We will obtain thus that $\tau$ is a symmetrising trace on
$\mathcal{H}(G_n)$. 
Note that we always have $1\in\mathcal{B}_n$, and that there exists a section $G_n \rightarrow \boldsymbol{G}_n$, $w \mapsto \boldsymbol{w}$ of $G_n$ in $B(G_n)$ such that  $\mathcal{B}_n = \{T_{\boldsymbol{w}}\,|\,\boldsymbol{w} \in \boldsymbol{G}_n\}$.
Therefore, $\tau$ satisfies Condition (1) of Conjecture \ref{BMM sym}.

By  definition, the linear map $\tau$ sends any element $h \in  \mathcal{H}(G_n)$ to the coefficient of $1$ when $h$ is expressed as a $\Z[\textbf{{u}},\textbf{{u}}^{-1}]$-linear combination of the elements of $\mathcal{B}_n$. We created a heuristic computer algorithm for a program in the language \texttt{C++} whose purpose was to produce this linear combination for every product $bb'$, with $b,b'\in \mathcal{B}_n$, and hence all entries of the matrix $A$.   However, with the exception of the case of $G_4$, this program was very time-consuming. So we created a second algorithm for a program in the language SAGE \cite{sagemath}   in order to optimise the calculation of $A$. The latter can be found on the project's webpage \cite{Chaw}, in a separate file for each group.

\subsection{The \texttt{C++} algorithm}\label{C++} 
The inputs of the \texttt{C++} algorithm are the following:  
 \begin{itemize}
\item[I1.] The basis $\mathcal{B}_n$.\smallbreak 
\item[I2.] The braid, positive and inverse Hecke relations  (cf.~Example \ref{Example of G4}). \smallbreak 
\item[I3.] The ``special cases'': these are some equalities computed by hand which express a given element of $\mathcal{H}(G_n)$ as a sum of other elements in $\mathcal{H}(G_n)$. 
In the appendix of this paper, we give the list of special cases for the group $G_6$, along with their proofs. For $G_4$ and $G_8$, the special cases are relatively easy, while for the  other groups, they are long and more complicated. The reader may find the lists of special cases for each group on the project's webpage \cite{Chaw}.
 \smallbreak 
\end{itemize}

Now, let $\mathcal{B}_n=\{b_1,\ldots,b_r\}$, and let $b,b' \in \mathcal{B}_n$. 
We know that the element $m:=bb'$ can be expressed uniquely as a linear combination $\lambda^m_1b_1 + \lambda^m_2 b_2 + \dots +\lambda^m_rb_r$ with $\lambda^m_i \in \Z[\textbf{{u}},\textbf{{u}}^{-1}]$. The algorithm takes the element $m$ and returns a list with the elements $\lambda^m_i b_i$ for which $\lambda^m_i  \neq 0$.

To start with, we define a list $L$ of  elements of the form $\lambda \,T_{\beta}$, where $\lambda \in \Z[\textbf{{u}},\textbf{{u}}^{-1}]$ and $\beta \in B(G_n)$. We  initialise $L$ with the element $m$, e.g. $L = [m]$. The algorithm described here is iterative. At each iteration, the initial element is  expanded and replaced by its summands.  For example, suppose that in $G_4$ the input is $m = s^2 \cdot s = s^3$. So in the beginning $L = [s^3]$. During the first iteration, by the positive Hecke relation, $s^3$ is equal to $as^2 + bs + c$, so the list $L$ is updated to become $L = [as^2, bs, c]$.

The same procedure is repeated at each iteration until all the elements of $L$ are of the form $\lambda \,T_{\beta}$ with $T_{\beta} \in \mathcal{B}_n$.
 When this is the case, the list $L$ is processed to sum up the coefficients for the same $T_{\beta}$, that is,
 if  $\lambda \,T_{\beta}, \lambda' \,T_{\beta} \in L$ for some $\beta \in B(G_n)$, then the two elements are removed from $L$ and
replaced by $(\lambda + \lambda')\,T_{\beta}$ only if $\lambda + \lambda' \neq 0$. For example, if in the end the list $L$ is $[(ab), (a^2c)st, (-a^2c)st, (bc), (bc^2)s]$, it will be processed for ``cleaning" the coefficients and become $L = [(ab + bc), (bc^2)s]$. 

We describe now what happens inside each iteration; this is the core of the algorithm.
For every element $l \in L$:
 \begin{itemize}
 	\item[S1.] We check whether $l\in \mathcal{B}_n$. If this is the case, there is nothing to do.	 \smallbreak
 	\item[S2.] If $l\not \in \mathcal{B}_n$, we check  whether $l$ appears in one of the special cases. If this happens, we replace inside $L$ the element $l$  by its summands given by this special case. For the groups $G_4$, $G_6$ and $G_8$, where only positive powers of the generators appear in the elements of $\mathcal{B}_n$, we use the inverse Hecke relations to transform any negative powers appearing in the special case into positive ones before replacing $l$. \smallbreak
 	\item[S3.] If $l\not  \in \mathcal{B}_n$ and it does not appear in any special case, we check whether we can apply the positive Hecke relation. If this is the case, we replace inside $L$ the element $l$ with its summands arising from the positive Hecke relation. \smallbreak
 	\item[S4.] If we are not able to apply any of the above steps, we create a second list $L'$ that contains  all elements that can transform into $l$ by applying the braid relation. We now check whether there is an element $l'$ in this list such that $l'\in \mathcal{B}_n$. If this is the case, we replace inside $L$ the element $l$ by the element $l'$. If not, we check whether there is an element $l' \in L'$ that appears in a special case. If the answer is positive, we replace again $l$ by $l'$ inside $L$ and we continue as in Step 2.
 \end{itemize}

\begin{rem}\rm The described steps correspond to the general form of the algorithm. Some adjustments, as for example the order of the steps, may be made depending on the group. In particular, Step 4  can be omitted if an important number of special cases has been developed, rendering this step unnecessary. 
\end{rem}

\subsubsection{\textbf{The group} $G_4$} As mentioned in Section \ref{Conjectures}, the BMM symmetrising trace conjecture has already been  proved for the group $G_4$, which  
is the smallest exceptional complex reflection group, in many different ways  \cite{MM10, Mar46}. The presentation of the 
generic Hecke algebra $\mathcal{H}(G_4)$ of $G_4$, together with the braid, positive and inverse Hecke relations, are given in Example \ref{Example of G4}.

Hoping that the appearance of positive powers of the generators $s$ and $t$ in the elements of $\mathcal{B}_4$ would make calculations easier, we changed the basis for $\mathcal{H}(G_4)$ given in 
\cite{Ch17} to the following one (the change is very straightforward, but it can be also justified by Theorem  \ref{Bn basis} in the next section):
$$\mathcal{B}_4=\{1,s,s^2,z,zs,zs^2,t^2,t,t^2s,ts,t^2s^2,ts^2,st^2,st,st^2s,sts,st^2s^2,sts^2,s^2t^2,s^2t,s^2t^2s,s^2ts,s^2t^2s^2,s^2ts^2\},$$
where $z:=T_{\boldsymbol{\rm z}}=ststst \in Z(\mathcal{H}(G_4))$. We have $|\mathcal{B}_4|=|G_4|=24$.

Using our program in \texttt{C++}, we were able to write every product $(bb')_{b,b'\in \mathcal{B}_4}$ as a $\mathbb{Z}[a,b,c^{\pm 1}]$-linear combination of elements in $\mathcal{B}_4$. The calculation of these linear combinations and, hence, of the matrix $A$ takes about 1 hour on an Intel Core i5 CPU.

\subsection{The SAGE algorithm}\label{SAGE}
Let $n \in \{5,\ldots,8\}$.
As we will see, the two algebra generators of $\mathcal{H}(G_n)$ are always in $\mathcal{B}_n$. Let us denote them by $g_1$ and $g_2$. We also have $z^k \in \mathcal{B}_n$, where $z:=T_{\boldsymbol{\rm z}} \in Z(\mathcal{H}(G_n))$ and $k=0,1,\ldots,|Z(G_n)|-1$. Note that
$z^{|Z(G_n)|}=T_{\boldsymbol{\pi}}$.

 The inputs of the SAGE algorithm  are the coefficients of the following elements when written as linear combinations of the elements of $\mathcal{B}_n$:
 \begin{itemize}
\item[I1.]  $g_1b$  for all $b \in \mathcal{B}_n$. \smallbreak
\item[I2.]  $g_2b$  for all $b \in \mathcal{B}_n$. \smallbreak
 \item[I3.] $z^{|Z(G_n)|} $.
  \end{itemize}
 The algorithm  then provides the matrix $A$ by  using the inductive nature of the basis $\mathcal{B}_n$. We observe that all inputs are computed by the \texttt{C++} algorithm, since all elements in the above list are products of two elements of the basis $\mathcal{B}_n$ (\emph{e.g.}, $z^{|Z(G_n)|}=z \cdot z^{|Z(G_n)|-1}$).
 
\begin{rem}\rm
For $n \in \{7,8\}$, input I3 is not even necessary, because it is deduced from the first two inputs (see Equations \eqref{eq3} and \eqref{z12} later on).
\end{rem}

\subsubsection{\textbf{The group} $G_8$}
We begin with the case of $G_8$, because it has the same braid relation as $G_4$.
The generic Hecke algebra of $G_8$ is the $\Z[u_0^{\pm 1},u_1^{\pm 1},u_2^{\pm 1},u_3^{\pm 1}]$-algebra:
$$\mathcal{H}(G_8)=\left\langle s,\,t \,\,\,\left|\,\,\, sts=tst,\,\,\, \prod_{i=0}^3(s-u_i) =  \prod_{i=0}^3(t-u_i) =0 \right\rangle\right..$$
We can rewrite the presentation of $\mathcal{H}(G_8)$ as follows:
$$\mathcal{H}(G_8)=\left\langle s,\,t \,\,\,\left|\,\,\, sts=tst,\,\,\, s^4 = as^3+bs^2+cs+d,\,\,\, t^4 = at^3+bt^2+ct+d\right\rangle\right..$$
Note that, in the above presentation, only $d$ is a unit in the ring of definition of $\mathcal{H}(G_8)$.
As in the case of $G_4$, besides the braid and the positive Hecke relations, we have the inverse Hecke relations 
\begin{center}
$s^{-1} =d^{-1}s^3-ad^{-1}s^2-bd^{-1}s-cd^{-1}$\,\,\, and \,\,\, $t^{-1} =d^{-1}t^3-ad^{-1}t^2-bd^{-1}t-cd^{-1}$.
\end{center} 

The basis for $\mathcal{H}(G_8)$ given in \cite{Ch17} is the following:
$$\mathcal{B}_8 =\left\{
\begin{matrix}
\begin{array}{c|l}
z^k, z^ks, z^ks^2, z^ks^3, z^kt, z^kts, z^kts^2, z^kts^3, z^kst, z^ksts, z^ksts^2, z^ksts^3, 
z^ks^2t,\;\;&\\
 z^ks^2ts, z^ks^2ts^2, z^ks^2ts^3, z^ks^3t, z^ks^3ts, z^ks^3ts^2, z^ks^3ts^3, 
z^kt^2, z^kst^2, z^ks^2t^2, z^ks^3t^2\;\;& 
\end{array}
k=0,1,2,3\;\;
\end{matrix}\right\},
$$
where $z=ststst \in Z(\mathcal{H}(G_8))$. We have $|\mathcal{B}_8|=|G_8|=96$ and $|Z(G_8)|=4$.
We now order the elements of $\mathcal{B}_8$. First,  we set
$$\mathcal{E}_8: =\left\{
\begin{matrix}
1, s, s^2, s^3, t, ts, ts^2, ts^3, st, sts, sts^2, sts^3, 
s^2t,
s^2ts, \\s^2ts^2, s^2ts^3, s^3t, s^3ts, s^3ts^2, s^3ts^3, 
t^2, st^2, s^2t^2, s^3t^2\;\;
\end{matrix}\right\}
\subset \mathcal{B}_8.$$
Every $i\in\{1,\dots, 96\}$ can be written as $24k+m$, where $k\in\{0,1,2,3\}$ and $m\in\{1,\dots, 24\}$. We write $\mathcal{B}_8=\{b_1,\dots,b_{96}\}$, where
$b_i:=z^kb_m$ for $b_m \in \mathcal{E}_8$ in the order written above. For example, $b_{72}=b_{24\cdot2+24}=z^2b_{24}=z^2s^3t^2$.
Note that for $1\leq i\leq 24$ we have $b_i\in \mathcal{E}_8$. 
We now notice that, for every $k\in\{0,1,2,3\}$, we have:
\begin{equation}
\label{r1}
\begin{array}{lll}
b_{24k+2\phantom{1}}= b_{24k+1\phantom{1}}\cdot s,\;\;\;&
b_{24k+3\phantom{1}}=  b_{24k+2\phantom{1}}\cdot s,&
\;\;\;b_{24k+4\phantom{1}}=  b_{24k+3\phantom{1}}\cdot s,\smallbreak\smallbreak\\
b_{24k+5\phantom{1}}= b_{24k+1\phantom{1}}\cdot t,&
b_{24k+6\phantom{1}}= b_{24k+5\phantom{1}}\cdot s,&
\;\;\;b_{24k+7\phantom{1}}= b_{24k+6\phantom{1}}\cdot s,\smallbreak\smallbreak\\
b_{24k+8\phantom{1}}= b_{24k+7\phantom{1}}\cdot s,&
b_{24k+9\phantom{1}}= b_{24k+2\phantom{1}}\cdot t,&
\;\;\;b_{24k+10}= b_{24k+9\phantom{1}}\cdot s,\smallbreak\smallbreak\\
b_{24k+11}= b_{24k+10}\cdot s,&
b_{24k+12}= b_{24k+11}\cdot s,&
\;\;\;b_{24k+13}= b_{24k+3\phantom{1}}\cdot t,\smallbreak\smallbreak\\
b_{24k+14}= b_{24k+13}\cdot s,&
b_{24k+15}= b_{24k+14}\cdot s,&
\;\;\;b_{24k+16}= b_{24k+15}\cdot s\smallbreak\smallbreak\\
b_{24k+17}= b_{24k+4\phantom{1}}\cdot t,&
b_{24k+18}= b_{24k+17}\cdot s,&
\;\;\;b_{24k+19}= b_{24k+18}\cdot s,\smallbreak\smallbreak\\
b_{24k+20}=  b_{24k+19}\cdot s,&
b_{24k+21}= b_{24k+5\phantom{1}}\cdot t,&
\;\;\;b_{24k+22}= b_{24k+9\phantom{1}}\cdot t,\smallbreak\smallbreak
\\
b_{24k+23}= b_{24k+13}\cdot t,&
b_{24k+24}=b_{12k+17}\cdot t.&
\end{array}
\end{equation}

Using the \texttt{C++} program, we have expressed 
$sb_{j}$ and $tb_{j}$ as $\mathbb{Z}[a,b,c,d^{\pm 1}]$-linear combinations of the elements of $\mathcal{B}_8$, for all $j\in\{1,\dots,96\}$. 
Let $sb_{j}:=\sum_{l}\lambda_{j,l}^{s}b_{l}$ and $tb_{j}:=\sum_{l}\lambda_{j,l}^{t}b_{l}$ be these linear combinations.  Following Equations \eqref{r1},  we can write every $b_{24k+m}$, with $m\in\{2,\dots,24\}$, as $b_{24k+m'}\cdot g$, with $ g \in\{s,t\}$ and $m'<m$.
Since, by definition, $\tau$ is  $\mathbb{Z}[a,b,c,d^{\pm 1}]$-linear, we obtain:
\begin{equation}\label{eq1}
\tau(b_{24k+m}b_{j})=\sum_{l} \lambda_{j,l}^{g} \tau(b_{24k+m'}b_{l}), \quad \text{ for all } j\in\{1,\dots,96\}.
\end{equation}

We now consider the case of $b_{24k+1}=z^k$, for $k\not=0$.
We distinguish two cases:
\begin{itemize}
	\item If $1 \leq j \leq 24(4-k)$, then $b_{24k+1}b_j=z^kb_j$.  Since $24k+1\leq 24k+j\leq 96$, we have $z^kb_j=b_{24k+j}$. We have assumed that $k\not=0$, so $b_{24k+j}\in \mathcal{B}_8 \setminus \{1\}$.
	Hence, 
	\begin{equation}\label{eq2}
	\tau(b_{24k+1}b_j )=0.
	\end{equation}
	\item If $24(4-k) < j \leq 96$, then $b_{24k+1}b_j=z^{k}b_j=
	z^{k-4} b_{j}z^4 = b_{24k+j-96}\cdot z^4$. Now notice that 
	$z^4= z^3tststs=tz^3s(tst)s= tz^3s^2ts^2=tb_{87}$.
	Hence, we already know how $z^4$ is written as a linear combination of the elements of $\mathcal{B}_8$, namely:
	\begin{equation}\label{eq3}
	z^4=\sum_{l}\lambda_{87,l}^{t}b_{l}.
	\end{equation}
	As a result, we obtain:
	\begin{equation}\label{eq4}
	\tau(b_{24k+1}b_j)=\tau(b_{24k+j-96}\cdot z^4)=\sum_{l}\lambda^{t}_{87,l}
	\tau(b_{24k+j-96}\,b_{l}).
	\end{equation} 
\end{itemize}

We can now compute the matrix $A$ row by row. 
Thanks to Equation \eqref{eq1}, with $m'$ and $g$ given by $\eqref{r1}$, one can fill in consecutively the entries of the rows $12k+2, 12k+3,\dots, 12k+24$, for every $k\in\{0,1,2,3\}$, with the only condition that the row $12k+1$ is known. Considering now the row $12k+1$ we notice the following:
The case $k=0$ corresponds to the first row of the matrix, which is given by $\tau(b_1b_j)=\delta_{1j}$, for all $j\in\{1,\dots, 96\}$.
For $k\not=0$ we use Equations \eqref{eq2} and \eqref{eq4}. Note that every term of the sum in \eqref{eq4} corresponds to rows 1 to $24k$, since $97-24k\leq j\leq 96$.

Summing up, the procedure is the following:
Knowing row 1, we can fill in consecutively the entries of rows 2--24. Since we know now rows 1--24, we can fill in row 25. Knowing now row 25 we can fill in 
consecutively the rows 26--48. From the rows 1--48 we can now fill in the entries of row 49. Knowing row 49 we can fill in consecutively the entries of the rows 50--72. With the entries of rows 1--72 we can fill in row 73 and, finally, we use row 73 to fill in consecutively the entries of the rows 74--96.

\subsubsection{\textbf{The group} $G_6$} The generic Hecke algebra of $G_6$ is the $\Z[u_{s,0}^{\pm 1},u_{s,1}^{\pm 1}, u_{t,0}^{\pm 1},u_{t,1}^{\pm 1},u_{t,2}^{\pm 1}]$-algebra:
$$\mathcal{H}(G_6)=\left\langle s,\,t \,\,\,\left|\,\,\, ststst=tststs,\,\,\, \prod_{i=0}^1(s-u_{s,i}) =  \prod_{j=0}^2(t-u_{t,j}) =0 \right\rangle\right..$$
We can rewrite the presentation of $\mathcal{H}(G_6)$ as follows:
$$\mathcal{H}(G_6)=\left\langle s,\,t \,\,\,\left|\,\,\, ststst=tststs,\,\,\, s^2 = as+b,\,\,\, t^3 = ct^2+dt+e\right\rangle\right..$$
Note that, in the above presentation, $b$ and $e$ are units in the ring of definition of $\mathcal{H}(G_6)$. 
As in the case of $G_8$, we have the braid relation, the positive Hecke relations and the inverse Hecke relations 
\begin{center}
$s^{-1}=b^{-1}s-ab^{-1}$ \,\,\, and  \,\,\, $t^{-1}=e^{-1}t^{2}-ce^{-1}t-de^{-1}$.
\end{center}

The basis 
 for $\mathcal{H}(G_6)$ given in \cite{Ch17}  is the following:
 $$\mathcal{B}_6 =\left\{
\begin{matrix}
z^k, z^kt, z^kt^2, z^ks, z^kst, z^kst^2, z^kts, z^kt^2s, z^ktst, z^ktst^2, z^kt^2st, z^kt^2st^2
\end{matrix}\;\;\left|\,\,k=0,1,2,3\right.\right\},
$$
where $z=ststst \in Z(\mathcal{H}(G_6))$.
We have $|\mathcal{B}_6|=|G_6|=48$ and $|Z(G_6)|=4$. Due to the inductive nature of the basis $\mathcal{B}_6$, we can again apply an algorithm similar to the one we used in the case of $G_8$. We start by ordering the elements of $\mathcal{B}_6$ as we did for the case of $G_8$. We set
$$\mathcal{E}_6: =\left\{
\begin{matrix}
1, t, t^2, s, st, st^2, ts, t^2s, tst, tst^2, t^2st, t^2st^2
\end{matrix}\right\}
\subset \mathcal{B}_6.$$
Every $i\in\{1,\dots, 48\}$ can be written as $12k+m$, where $k\in\{0,1,2,3\}$ and $m\in\{1,\dots, 12\}$. We write $\mathcal{B}_6=\{b_1,\dots,b_{48}\}$, where
$b_i:=z^kb_m$ for $b_m \in \mathcal{E}_6$ in the order written above.
We now notice that, for every $k\in\{0,1,2,3\}$, we have:
\begin{equation}
\label{r11}
\begin{array}{lll}
b_{12k+2\phantom{1}}=b_{12k+1\phantom{1}}\cdot t,\;\;&
b_{12k+3\phantom{1}}=  b_{12k+2\phantom{1}}\cdot t,\;\;&
b_{12k+4\phantom{1}}=b_{12k+1\phantom{1}}\cdot s,\smallbreak\smallbreak\\
b_{12k+5\phantom{1}}=b_{12k+4\phantom{1}}\cdot t,&
b_{12k+6\phantom{1}}= b_{12k+5\phantom{1}}\cdot t,&
b_{12k+7\phantom{1}}= b_{12k+2\phantom{1}}\cdot s,\smallbreak\smallbreak\\
b_{12k+8\phantom{1}}= b_{12k+3\phantom{1}}\cdot s,&
b_{12k+9\phantom{1}}= b_{12k+7\phantom{1}}\cdot t,&
b_{12k+10}= b_{12k+9\phantom{1}}\cdot t,\smallbreak\smallbreak\\
b_{12k+11}= b_{12k+8\phantom{1}}\cdot t,&
b_{12k+12}= b_{12k+11}\cdot t.&
\end{array}
\end{equation}

Using the \texttt{C++} program we have expressed, for all $j\in\{1,\dots,48\}$, 
$sb_{j}$, $tb_{j}$ and $z^4=b_{13}b_{37}$ as linear combinations of the elements of $\mathcal{B}_6$ with coefficients in $\mathbb{Z}[a,b^{\pm1},c,d,e^{\pm1}]$.
Let $sb_{j}:=\sum_{l}\lambda^s_{j,l}b_{l}$,  $tb_{k,j}:=\sum_{l}\lambda^t_{j,l}b_{l}$ and $z^4:=\sum_{l}\mu_{l}b_{l}$
be these linear combinations. 
Thanks to Equations \eqref{r11}, we can write every $b_{12k+m}$, with $m\in\{2,\dots,12\}$, as $b_{12k+m'}\cdot g$, where $g\in\{s,t\}$ and $m'<m$. Since, by definition, $\tau$ is  $\mathbb{Z}[a,b^{\pm1},c,d,e^{\pm1}]$-linear, we obtain:
$$
\tau(b_{12k+m}b_{j})=\sum_{l} \lambda^{g}_{j,l} \tau(b_{12k+m'}b_{l}),
\quad \text{ for all } j \in \{1,\ldots,48\}.
$$

We now consider the case of $b_{12k+1}=z^k$, for $k\not=0$.
We distinguish two cases, as we did for $G_8$:
\begin{itemize}
	\item If $1 \leq j \leq 12(4-k)$, then we have
	$\tau(b_{24k+1}b_j )=0$. \smallbreak
	\item If $12(4-k) < j \leq 48$, then $b_{12k+1}b_j=z^{k}b_j=
	z^{k-4} b_{j}z^4 = b_{12k+j-48}\cdot z^4$. 
	As a result, we have:
	$$
	\tau(b_{12k+1}b_j)=\tau(b_{12k+j-48}\cdot z^4)=\sum_{l}\mu_l
	\tau(b_{12k+j-48}\,b_{l}).$$
\end{itemize}
We fill in the entries of the matrix $A$ row by row, as we did in the case of $G_8$.

\subsubsection{\textbf{The group} $G_5$}\label{G5} The generic Hecke algebra of $G_5$ is the $\Z[u_{s,0}^{\pm 1},u_{s,1}^{\pm 1},u_{s,2}^{\pm 1}, u_{t,0}^{\pm 1},u_{t,1}^{\pm 1},u_{t,2}^{\pm 1}]$-algebra:
$$\mathcal{H}(G_5)=\left\langle s,\,t \,\,\,\left|\,\,\, stst=tsts,\,\,\, \prod_{i=0}^2(s-u_{s,i}) =  \prod_{j=0}^2(t-u_{t,j}) =0 \right\rangle\right..$$
We can rewrite the presentation of $\mathcal{H}(G_5)$ as follows:
$$\mathcal{H}(G_5)=\left\langle s,\,t \,\,\,\left|\,\,\, stst=tsts,\,\,\, s^3 = as^2+bs+c,\,\,\, t^3 = dt^2+et+f\right\rangle\right..$$
Note that, in the above presentation, $c$ and $f$ are units in the ring of definition of $\mathcal{H}(G_5)$. We have again the braid relation, the positive Hecke relations, and the inverse Hecke relations
\begin{center}
 $s^{-1}=c^{-1}t^2-ac^{-1}t-bc^{-1}$ \,\,\, and \,\,\, $t^{-1}=f^{-1}t^2-df^{-1}t-ef^{-1}$.
 \end{center}
 
The basis 
 for $\mathcal{H}(G_5)$ given in \cite{Ch17}  is the following one:
 $$\mathcal{B}_5 =\left\{
\begin{matrix}
z^k, z^ks, z^ks^2, z^kt, z^kt^2, z^kst,  z^ks^2t, z^kst^2, z^ks^2t^2, z^kt^{-1}s, z^kt^{-1}st, z^kt^{-1}st^2
\end{matrix}\,\,\left|\,\,k=0,1,2,3,4,5\right.\right\},
$$
where $z=stst \in Z(\mathcal{H}(G_5))$.
We have $|\mathcal{B}_5|=|G_5|=72$ and $|Z(G_5)|=6$. Due to the inductive nature of the basis $\mathcal{B}_5$, we can again apply an algorithm similar to the one we used in the cases of $G_8$ and $G_6$.
We order the elements of $\mathcal{B}_5$ as in the previous cases, by setting
$$\mathcal{E}_5: =\left\{
\begin{matrix}
1, s, s^2, t, t^2, st,  s^2t, st^2, s^2t^2, t^{-1}s, t^{-1}st, t^{-1}st^2
\end{matrix}\right\}
\subset \mathcal{B}_5.$$
Every $i\in\{1,\dots, 72\}$ can be written as $12k+m$, where $k\in\{0,\ldots,5\}$ and $m\in\{1,\dots, 12\}$. We write $\mathcal{B}_5=\{b_1,\dots,b_{72}\}$, where
$b_i:=z^kb_m$ for $b_m \in \mathcal{E}_5$ in the order written above.

We will now provide a set of inductive relations, analogous to \eqref{r1} and \eqref{r11}.
However, the presence of a negative power $t^{-1}$ in the elements $b_{12k+10}$, $k\in\{0,\dots,5\}$, complicates the 10th relation, compared to the others which are straightforward (note that this negative power does not affect the relations for the elements $b_{12k+11}$ and $b_{12k+12}$). 
 We use the inverse Hecke relation that corresponds to $t$ and we obtain:
 $$b_{12k+10}=z^kb_{10}=z^kt^{-1}s=f^{-1}(z^kt^2-dz^kt-ez^k)s.$$ As a result, for every $k\in\{0,\dots,5\}$, we have:
\begin{equation}
\label{r111}
\begin{array}{lll}
b_{12k+2\phantom{1}}= b_{12k+1\phantom{1}}\cdot s,&
b_{12k+3\phantom{1}}=  b_{12k+2\phantom{1}}\cdot s,&
b_{12k+4\phantom{1}}=  b_{12k+1\phantom{1}}\cdot t,\smallbreak\smallbreak\\
b_{12k+5\phantom{1}}= b_{12k+4\phantom{1}}\cdot t,&
b_{12k+6\phantom{1}}= b_{12k+2\phantom{1}}\cdot t,&
b_{12k+7\phantom{1}}=b_{12k+3\phantom{1}}\cdot t,\smallbreak\smallbreak\\
b_{12k+8\phantom{1}}=b_{12k+6\phantom{1}}\cdot t,&
b_{12k+9\phantom{1}}=b_{12k+7\phantom{1}}\cdot t,&
b_{12k+10}=f^{-1}(b_{12k+5}-db_{12k+4}-eb_{12k+1}) \cdot s,\smallbreak\smallbreak\\
b_{12k+11}= b_{12k+10}\cdot t,&
b_{12k+12}= b_{12k+11}\cdot t.&
\end{array}
\end{equation}

Using the \texttt{C++} program, we have expressed, for all $j\in\{1,\dots,72\}$,
$sb_{j}$, $tb_{j}$ and $z^6=
b_{37}^2$ as linear combinations of the elements of $\mathcal{B}_5$ with coefficients in $\mathbb{Z}[a,b,c^{\pm1},d,e,f^{\pm1}]$.
Let $sb_{j}:=\sum_{l}\lambda^s_{j,l}b_{l}$,  $tb_{j}:=\sum_{l}\lambda^t_{j,l}b_{l}$ and $z^6:=\sum_{l}\mu_{l}b_{l}$
be these linear combinations. 
Thanks  to Equations \eqref{r111}, we can write every $b_{12k+m}$, with $m\in\{2,\dots,12\}\setminus\{10\}$, as $b_{12k+m'}\cdot g$, where $g \in\{s,t\}$ and $m'<m$. 
 Since, by definition, $\tau$ is  $\mathbb{Z}[a,b,c^{\pm1},d,e,f^{\pm1}]$-linear, we obtain:
$$
\tau(b_{12k+m}b_{j})=\sum_{l} \lambda^{g}_{j,l} \tau(b_{12k+m'}b_{l}), \quad \text{ for all } j\in\{1,\dots,72\}.
$$
For $m=10$, we have: 
$$\tau (b_{12k+10}b_j)= f^{-1}\sum_{l}\lambda^s_{j,l}(\tau(b_{12k+5}b_l) -d\tau(b_{12k+4}b_l)-
e\tau(b_{12k+1}b_l)).$$

We now consider the case of $b_{12k+1}=z^k$, for $k\not=0$.
We distinguish two cases, as we did for $G_8$ and $G_6$:
\begin{itemize}
	\item If $1 \leq j \leq 12(6-k)$, then we have
	$\tau(b_{12k+1}b_j )=0$. \smallbreak
	\item If $12(6-k)< j \leq 72$, then $b_{12k+1}b_j=z^{k}b_j=
	z^{k-6} b_{j}z^6 = b_{12k+j-72}\cdot z^6$. 
	As a result, we have:
	$$
	\tau(b_{12k+1}b_j)=\tau(b_{12k+j-72}\cdot z^6)=\sum_{l}\mu_l
	\tau(b_{12k+j-72}\, b_{l}).$$
\end{itemize}
We fill in the entries of the matrix $A$ row by row, as we did for the cases of $G_8$ and $G_6$.

\subsubsection{\textbf{The group} $G_7$}\label{G7} The generic Hecke algebra of $G_7$ is the $\Z[u_{s,0}^{\pm 1},u_{s,1}^{\pm 1}, u_{t,0}^{\pm 1},u_{t,1}^{\pm 1},u_{t,2}^{\pm 1},u_{u,0}^{\pm 1},u_{u,1}^{\pm 1},u_{u,2}^{\pm 1}]$-algebra:
$$\mathcal{H}(G_7)=\left\langle s,\,t,\,u \,\,\,\left|\,\,\, stu=tus=ust,\,\,\, \prod_{i=0}^1(s-u_{s,i}) =  \prod_{j=0}^2(t-u_{t,j}) =\prod_{k=0}^2(u-u_{u,j})=0 \right\rangle\right..$$
We can rewrite the presentation of $\mathcal{H}(G_7)$ as follows:
$$\mathcal{H}(G_7)=\left\langle s,\,t,\,u \,\,\,\left|\,\,\, stu=tus=ust,\,\,\, s^2 = as+b,\,\,\, t^3 = ct^2+dt+e,\,\,\,u^3=fu^2+gu+h\right\rangle\right..$$
Note that, in the above presentation, $b$, $e$ and $h$ are units in the ring of definition of $\mathcal{H}(G_7)$. As in the previous cases, we have the braid relations, the positive Hecke relations, and the inverse Hecke relations 
\begin{center}
$s^{-1}=b^{-1}s-ab^{-1}$,\,\,\, $t^{-1} = e^{-1}t^2 - ce^{-1}t - de^{-1}$ \,\,\, and \,\,\, $u^{-1} =h^{-1}u^2 - fh^{-1}u - gh^{-1}$. 
\end{center}

The basis 
 for $\mathcal{H}(G_7)$ given in \cite{Ch17}  is the following one:
 $$\widetilde{\mathcal{B}}_7 =
 \left\{
\begin{matrix}
z^k, z^ku, z^ku^2, z^kt, z^kt^2, z^kut,  z^ku^2t, z^kut^2, z^ku^2t^2, z^ktu^{-1}, z^ktu^{-1}t, z^ktu^{-1}t^2
\end{matrix}\,\,\left|\,\,k=0,1,\ldots,11\right.\right\},
$$
where $z=tus \in Z(\mathcal{H}(G_7))$.
We have $|\widetilde{\mathcal{B}}_7 |=|G_7|=144$ and $|Z(G_7)|=12$. 
Unfortunately, we quickly discovered that by setting
$\tau(b):=\delta_{1b}$ for all $b \in \widetilde{\mathcal{B}}_7$ we do not obtain a trace function. For example, let us calculate the trace of the elements $u^2t^2\cdot t$ and
$t\cdot u^2t^2$. We have $\tau(u^2t^3)=c\tau(u^2t^2)+d\tau(u^2t)+e\tau(u^2)=0$. In order to calculate
$\tau(tu^2t^2)$, we first multiply the positive Hecke relation for $u$ by $u^{-1}$ to obtain the ``equivalent positive Hecke relation''
$u^2=fu+g+hu^{-1}$. We then have
$$\tau(tu^2t^2)=f\tau(tut^2)+g\tau(t^3)+h\tau(tu^{-1}t^2)=f\tau(tut^2)+g\tau(t^3).$$
Since $t^3=ct^2+dt+e$, we have $\tau(t^3)=e$. 
Now, $tut^2=(tus)s^{-1}t^2=zs^{-1}t^2$. Using the inverse Hecke relation for $s^{-1}$ yields 
$$zs^{-1}t^2=b^{-1}zst^2-ab^{-1}zt^2=b^{-1}z(stu)u^{-1}t-ab^{-1}zt^2=b^{-1}z^2u^{-1}t-ab^{-1}zt^2.$$ Using the inverse Hecke relation for $u^{-1}$ yields
$$zs^{-1}t^2=b^{-1}h^{-1}z^2u^{2}t - b^{-1}fh^{-1}z^2ut -b^{-1}gh^{-1}z^2t- ab^{-1}zt^2,$$
 whence we deduce that $\tau(tut^2)= \tau(zs^{-1}t^2)=0$. As a result, $\tau(tu^2t^2)=ge \neq 0 =\tau(u^2t^3)$.

One may suggest that the above problem could be solved by replacing the elements
$z^ktu^{-1}$, $z^ktu^{-1}t$, $z^ktu^{-1}t^2$ inside $\widetilde{\mathcal{B}}_7$ with
$z^ktu^{2}$, $z^ktu^{2}t$, $z^ktu^{2}t^2$ respectively, for all $k=0,1,\ldots,11$.
Then, by definition of $\tau$, we would have $\tau(tu^2t^2)=0$.
Unfortunately, the set obtained this way is not a basis for $\mathcal{H}(G_7)$, and this was the reason behind the introduction of negative powers of the generators in the basis given in \cite{Ch17}.

Nevertheless, it turns out that replacing only the three elements $tu^{-1}$, $tu^{-1}t$, $tu^{-1}t^2$ inside $\widetilde{\mathcal{B}}_7$ with
$tu^{2}$, $tu^{2}t$, $tu^{2}t^2$ respectively does yield a basis for $\mathcal{H}(G_7)$, and rectifies the problem caused by $tu^{2}t^2$. 
The reason we replace three elements and not only the element $tu^2t^2$ is to be able to apply the SAGE algorithm as we did in the previous cases. 
Therefore, we define
$$\mathcal{B}_7: =\left( \widetilde{\mathcal{B}}_7 \setminus \{tu^{-1}, tu^{-1}t, tu^{-1}t^2\} \right)
\cup \{tu^{2}, tu^{2}t, tu^{2}t^2\}.$$

The set $\mathcal{B}_7$ is a basis for $\mathcal{H}(G_7)$, because, for all $r=0,1,2$, we have:
\begin{equation}
\begin{array}{c}\smallbreak \smallbreak \smallbreak
tu^{-1}t^r = h^{-1}tu^2t^r - fh^{-1}tut^r - gh^{-1}t^{r+1} \\ \smallbreak \smallbreak \smallbreak
tut^r = s^{-1}(stu)t^r = s^{-1}zt^r = b^{-1}szt^r-ab^{-1}zt^r\\ \smallbreak \smallbreak \smallbreak
szt^r=zst^r=z(stu)u^{-1}t^{r-1}=z^2u^{-1}t^{r-1}=h^{-1}z^2u^2t^{r-1} - fh^{-1}z^2ut^{r-1} - gh^{-1}z^2t^{r-1}.
\end{array}
\end{equation}
Now, given the positive and inverse Hecke relations for $t$,  it is easy to see that
$tu^{-1}t^r$ can be expressed as a linear combination of elements of $\mathcal{B}_7$, for all $r=0,1,2$.  Another proof of the fact that $\mathcal{B}_7$ is a basis for $\mathcal{H}(G_7)$  is given in Section \ref{Results} (cf.~Theorem \ref{Bn basis}).
In case the reader is wondering what happens with the algebra generator $s$, we observe that $s$ is written as a linear combination of the elements of $\mathcal{B}_7$ as follows:
$$
s =(stu)u^{-1}t^{-1} = zu^{-1}t^{-1} = e^{-1}h^{-1} (zu^2t^2 - czu^2t - fzut^2 - dzu^2 - gzt^2 +cf zut  + dfzu + cgzt +dg z).
$$

We order the elements of $\mathcal{B}_7$  by setting
$$\mathcal{E}_7: =\left\{
\begin{matrix}
1, u, u^2, t, t^2, ut,  u^2t, ut^2, u^2t^2, tu^{-1}, tu^{-1}t, tu^{-1}t^2
\end{matrix}\right\}.$$
Every $i\in\{1,\dots, 144\}$ can be written as $12k+m$, where $k\in\{0,\dots,11\}$ and $m\in\{1,\dots, 12\}$. 
We write $\mathcal{B}_7=\{b_1,\dots,b_{144}\}$, where
\begin{itemize}
\item for every $i\in\{1,\dots, 144\}\setminus\{10,11,12\}$,
$b_i:=z^kb_m$ for $b_m \in \mathcal{E}_5$ in the order written above;\smallbreak  
\item $b_{10}:=tu^{2}$, $b_{11}:=tu^{2}t$, $b_{12}:=tu^{2}t^2$. 
\end{itemize}

We now provide again a set of inductive relations,  analogous to \eqref{r1}, \eqref{r11} and \eqref{r111}.
The only relation which is not as in the previous cases is the one that corresponds to the elements $b_{12k+10}$. More precisely, if $k\not=0$, then $b_{12k+10}=z^ktu^{-1}$.
We use the inverse Hecke relation that corresponds to $u$ and we obtain:
$$b_{12k+10}=z^ktu^{-1}=h^{-1}z^ktu^2-fh^{-1}z^ktu-gh^{-1}z^kt.$$ 
As a result, for every $k\in\{0,\dots,11\}$, we have:
\begin{equation}
\label{r1111}
\begin{array}{lll}
b_{12k+2\phantom{1}}=b_{12k+1\phantom{1}}\cdot u,&
b_{12k+3\phantom{1}}= b_{12k+2\phantom{1}}\cdot u,&
b_{12k+4\phantom{1}}=  b_{12k+1\phantom{1}}\cdot t,\smallbreak\smallbreak\\
b_{12k+5\phantom{1}}= b_{12k+4\phantom{1}}\cdot t,&
b_{12k+6\phantom{1}}= b_{12k+2\phantom{1}}\cdot t,&
b_{12k+7\phantom{1}}= b_{12k+3\phantom{1}}\cdot t,\smallbreak\smallbreak\\
b_{12k+8\phantom{1}}=b_{12k+6\phantom{1}}\cdot t,&
b_{12k+9\phantom{1}}=b_{12k+7\phantom{1}}\cdot t,&
b_{12k+10}=\begin{cases}
tu^2 & \text{if } k=0 \\
h^{-1}(z^ktu^2-fz^ktu-gz^kt) & \text{if } k\not=0
\end{cases}, \smallbreak\smallbreak\\
b_{12k+11}=b_{12k+10}\cdot t,&
b_{12k+12}=b_{12k+11}\cdot t.&
\end{array}
\end{equation}
The above relations do not seem inductive for $b_{12k+10}$. However, we will see in a while that we do have an inductive relation for $\tau(b_{12k+10}b_j)$, $j\in\{1,\dots,144\}$, which is what we really want in order to apply the SAGE algorithm.

Using the \texttt{C++} program, we have expressed 
$ub_{j}$ and $tb_{j}$ as $\mathbb{Z}[a,b^{\pm},c,d,e^{\pm1},f,g,h^{\pm1}]$-linear combinations of the elements of $\mathcal{B}_7$ for all $j\in\{1,\dots,144\}$. 
Let $ub_{j}:=\sum_{l}\lambda^u_{j,l}b_{l}$ and $tb_{j}:=\sum_{l}\lambda^t_{j,l}b_{l}$
be these linear combinations. 
The third input of the SAGE algorithm are the coefficients of the element $z^{12}$ when expressed as a linear combination of the elements of the basis. However, as in the case of $G_8$ (cf.~Equation \eqref{eq3}), we do not need the \texttt{C++} program in order to calculate these coefficients on top of the ones we have already determined, because we notice the following:
\begin{equation}\label{z12}
\begin{array}{lcl}\smallbreak\smallbreak
z^{12}&=&z^{11}ust
\,=\,z^{11}u(a+bs^{-1})t
\,=\,a\,z^{11}ut+b\,z^{11}u(s^{-1}u^{-1}t^{-1})tut 
\,=\,a\,b_{138}+b\, z^{10}utut\\  
&=&a\, b_{138}+b\,utz^{10}ut
\,=\,a\, b_{138}+b\, utb_{126}
\,=\,a\, b_{138}+b\, u\sum\limits_{r}\lambda^t_{126,r}b_{r}\\
&=&a\, b_{138}+b\,\sum\limits_{r}\lambda^t_{126,r}ub_{r}
\,=\,a\, b_{138}+b\,\sum\limits_{r}\lambda^t_{126,r}\sum\limits_{l}\lambda^u_{r,l}b_{l}\,.
	\end{array}
	\end{equation}
	Hence, we have $z^{12}=\sum_{l}\mu_{l}b_{l}$, where
$\mu_{l}:=
a\,\delta_{l,138}+b\,\sum\limits_{r}\lambda^t_{126,r}\sum\limits_{l}\lambda^u_{r,l}$.

Due to Equations \eqref{r1111} we can write every $b_{12k+m}$, with $m\in\{2,\dots,12\}\setminus\{10\}$, as $b_{12k+m'}\cdot g$, where $g \in\{u,t\}$ and $m'<m$. 
Since, by definition, $\tau$ is  $\mathbb{Z}[a,b^{\pm1},c,d,e^{\pm1},f,g,h^{\pm1}]$-linear, we obtain:
\begin{equation}
\tau(b_{12k+m}b_{j})=\sum_{l} \lambda^{g}_{j,l} \tau(b_{12k+m'}b_{l}), \quad \text{ for all } j\in\{1,\dots,144\}.
\label{t1}
\end{equation}
For $m=10$, we distinguish two cases: the case $k=0$ and the case $k\not=0$. In both cases we will need the following calculation: for all $j\in\{1,\dots,144\}$, we have
\begin{equation}\label{u2}
u^2b_j = u\sum\limits_{l}\lambda^u_{j,l}b_{l} = \sum\limits_{l}\lambda^u_{j,l}(ub_{l}) = \sum\limits_{l}\lambda^u_{j,l}\sum\limits_{r}\lambda^u_{l,r}b_r. 
\end{equation}
Now,
\begin{itemize}[leftmargin=*]
	\item if $k=0$, then, using \eqref{u2}, we obtain
		$$b_{10}b_j= tu^2b_j = \sum\limits_{l}\lambda^u_{j,l}
		\sum\limits_{r}\lambda^u_{l,r}(tb_r)=\sum\limits_{l}\lambda^u_{j,l}
	\sum\limits_{r}\lambda^u_{l,r}
	\sum\limits_{p}\lambda^t_{r,p}b_p,$$
	and so,
		\begin{equation}
			\label{t2}
				\tau(b_{10}b_j)=\sum\limits_{l,r}\lambda^u_{j,l}
	\lambda^u_{l,r}
	\lambda^t_{r,1}. 
\end{equation}
		\item if $k\not=0$, then, 
		by \eqref{r1111}, we have $b_{12k+10}
		=h^{-1}z^kt \,(u^2-fu-g)$. Hence, using \eqref{u2}, we obtain
		$$\begin{array}[t]{lcl}\smallbreak\smallbreak\smallbreak
		b_{12k+10}b_j&=&h^{-1}z^kt\,(u^2b_j-fub_j-gb_j)\\
		&=&h^{-1}b_{12k+4} \left(\sum\limits_{l}\lambda^u_{j,l}\sum\limits_{r}\lambda^u_{l,r}b_r 
-f \sum\limits_{l}\lambda^u_{j,l}b_l -gb_j\right)
				\smallbreak\smallbreak\smallbreak\\
		&=&h^{-1}\sum\limits_{l}\lambda^u_{j,l}
		\sum\limits_{r}\lambda^u_{l,r}b_{12k+4}b_r
		-fh^{-1}\sum\limits_{l}\lambda^u_{j,l}b_{12k+4}b_{l}
			-gh^{-1} b_{12k+4}b_j,
		\end{array}$$ and so,
		\begin{equation}
		\label{t3}
		\tau(b_{12k+10}b_j)=h^{-1}\sum\limits_{l,r}\lambda^u_{j,l}
	\lambda^u_{l,r} \tau(b_{12k+4}b_r)
		-fh^{-1}\sum\limits_{l}\lambda^u_{j,l}\tau(b_{12k+4}b_{l})
		-gh^{-1}\tau(b_{12k+4}b_j).
		\end{equation}
\end{itemize}

We now consider the case of $b_{12k+1}=z^k$, for $k\not=0$. Since our basis does not have the pattern we had in the previous cases, we need to do some extra work:
\begin{itemize}
	\item 
	If $1 \leq j \leq 12(12-k)$, then
	we have
	$b_{12k+1}b_j=z^kb_j$. For $j\not\in\{10,11,12\}$,  $z^kb_j\in\mathcal{B}_7\setminus\{1\}$ and hence, $\tau(b_{12k+1}b_j)=0$.
	However, for $j\in\{10,11,12\}$, we have $b_{12k+1}b_j=z^ktu^2t^r$, for $r=0,1,2$. We use the equivalent positive Hecke relation $u^2=fu+g+hu^{-1}$ and we have:
	$$z^ktu^2t^r=fz^ktut^r+gz^kt^{r+1}+hz^ktu^{-1}t^r.$$
We notice now that, since $k\not=0$,  $z^ktu^{-1}t^r\in \mathcal{B}_7\setminus\{1\}$. The same holds for $z^kt^{r+1}$, when $r\not=2$. For $r=2$ we have $z^kt^3=cz^kt^2+dz^kt+ez^k$ and, hence, since $k\not=0$, $z^kt^3$ is a linear combination of elements in 
	$\mathcal{B}_7\setminus\{1\}$. As a result, $\tau(z^ktu^{-1}t^r)=\tau(z^kt^{r+1})=0$.
	 We focus now on $z^ktut^r$. We have 
	 $$ z^ktut^r=tz^kut^r= \left\{ \begin{array}{ll}
	 tb_{12k+2} & \text{ if } r=0,\\
	 tb_{12k+6} & \text{ if } r=1,\\
	 tb_{12k+8} & \text{ if } r=2.\\
	 \end{array}\right.$$
	 	 According to the \texttt{C++} program, $\lambda^t_{12k+m,1}=0$, for every $k\not=0$ and $m\in\{2,6,8\}$. As a result,
	  $\tau(z^ktut^r)=0$  and hence, $\tau(z^ktu^2t^r)=0$, for all $r\in\{0,1,2\}$.
	 Summing up, we have:
\begin{equation}
\label{t4}
\tau(b_{12k+1}b_j )=0, \quad \text{ for all } 1 \leq j \leq 12(12-k).
\end{equation}
	\item 
	If $12(12-k) < j \leq 144$, then,
	as in all the previous cases of $G_8$, $G_6$ and $G_5$, we have $b_{12k+1}b_j=z^{k}b_j=
	z^{k-12} b_{j}z^{12}$. However, the element $z^{k-12}b_{j}$ does not correspond every time to an element in $\mathcal{B}_7$, as before. This problem occurs when
	$z^{k-12}b_j=tu^{-1}t^r$ for $r=0,1,2$, that is, when $j=154-12k, 155-12k, 156-12k$ respectively.
	Therefore, if $j$ does not take one of these three values, then
	\begin{equation}
		\label{t5}
		\tau(b_{12k+1}b_j)=
		\tau(b_{12k+j-144}\cdot z^{12}) = \sum_l \mu_l \tau(b_{12k+j-144}b_l).
				\end{equation}
	If now $j \in  \{154-12k, 155-12k, 156-12k\}$, then, using the inverse Hecke relation for $u$, we obtain 
	$$
		b_{12k+1}b_j=
		tu^{-1}t^r z^{12} = h^{-1}tu^2t^rz^{12}-h^{-1}ftut^rz^{12}-h^{-1}gt^{r+1}z^{12}
		\quad \text{ for } r=0,1,2.
				$$
	The latter is equal to
	$$h^{-1} \sum_l \mu_l \,(tu^2t^rb_l-ftut^rb_l-gt^{r+1}b_l).$$
	We now have
	$$ tu^2t^rb_l= \left\{ \begin{array}{ll}
	 b_{10}b_l & \text{ if } r=0,\\
	 b_{11}b_l & \text{ if } r=1,\\
	 b_{12}b_l & \text{ if } r=2,\\
	 \end{array}\right.
	 \quad
	 tut^rb_l= \left\{ \begin{array}{ll}
	 \sum_{p} \lambda_{l,p}^u b_4b_p & \text{ if } r=0,\\
	 \sum_{p,q} \lambda_{l,p}^t\lambda_{p,q}^u b_4b_q & \text{ if } r=1,\\
	 \sum_{p,q,x} \lambda_{l,p}^t\lambda_{p,q}^t\lambda_{q,x}^u b_4b_x& \text{ if } r=2,\\
	 \end{array}\right.
	$$ $$
	 t^{r+1}b_l= \left\{ \begin{array}{ll}
	 b_{4}b_l & \text{ if } r=0,\\
	 b_{5}b_l & \text{ if } r=1,\\
	 cb_{5}b_l + db_4b_l +eb_l & \text{ if } r=2,\\
	 \end{array}\right.
	 $$
	for all $l \in \{1,\ldots,144\}$. We conclude that:\\
	$\star$ If $j=154-12k$, then $\tau(b_{12k+1}b_j)=
		\tau(tu^{-1}z^{12})$, which in turn is equal to
	\begin{equation}
	\label{t6} 	 
	h^{-1}\sum\limits_{l}\mu_{l} \left( \tau(b_{10}b_l)
	-f\sum\limits_{p}\lambda^u_{l,p}\tau(b_4b_p)
	-g\tau(b_4b_l)\right).
	\end{equation}
		$\star$ If $j=155-12k$, then $\tau(b_{12k+1}b_j)=
		\tau(tu^{-1}tz^{12})$, which in turn is equal to
		\begin{equation}
		\label{t7}
		h^{-1}\sum\limits_{l}\mu_{l} \left(\tau(b_{11}b_l)
		-f\sum\limits_{p,q}\lambda^t_{l,p}
		\lambda^u_{p,q}\tau(b_4b_q)
		-g\tau(b_5b_l)\right).
		\end{equation}
			$\star$ If $j=156-12k$, then $\tau(b_{12k+1}b_j)=
			\tau(tu^{-1}t^2z^{12})$, which in turn is equal to
\begin{equation}
		\label{t8}
h^{-1}\sum\limits_{l}\mu_{l}
\left( \tau(b_{12}b_l)
			- f\sum\limits_{p,q,x}\lambda_{l,p}^t\lambda_{p,q}^t\lambda_{q,x}^u\tau(b_4b_x)
			-gc \tau(b_5b_l)-gd\tau(b_4b_l)-ge\tau(b_l)\right).\end{equation}
			\end{itemize}

Let us now compute the matrix $A$. 
From \eqref{r1111}, \eqref{t1}, \eqref{t2}, and \eqref{t3}
we can complete consecutively rows $12k+2$ to $12k+12$ if we know row $12k+1$.
Considering now row $12k+1$, we notice that the first row ($k=0$) is known, since $\tau(b_1b_j)=\delta_{1j}$. For $k\not=0$, we can use 
\eqref{t4}, \eqref{t5}, \eqref{t6}, \eqref{t7} and \eqref{t8}, since every term of the sums in these equations correspond to rows 1 to $12k$.
Hence, we can fill in the entries of the matrix $A$ row by row, as we did for the cases of $G_5$, $G_6$, and $G_8$.

\subsection{The  GAP3 algorithm}\label{GAP3} Here we present another algorithm for the determination of the matrix $A$, which is easy to come up with but turned out to be computationally inefficient for the three larger groups. The program based on this algorithm produced the matrix $A$ for the groups $G_4$ and $G_6$, but it did not even compute the first line for the groups $G_5$ and $G_7$.
As far as $G_8$ is concerned, it produced the matrix $A$, but it could not simplify the entries enough to allow the calculation of the determinant.
Nevertheless, the program worked for all groups  when we replaced the parameters $a$, $b$, $c$, etc. with distinct prime numbers, and its results for the matrix $A$ and the determinant coincided with the ones we got from the \texttt{C++}/SAGE program. 

The algorithm makes use of the data on the Hecke algebras of complex reflection groups that are stored in the GAP3 package CHEVIE \cite{chevie1,chevie2}. It consists of GAP3 functions, but, in the end, we ran it on \emph{Mathematica}, because it is faster and more efficient with large polynomials than GAP3.

 Now, in  CHEVIE, we can find the Schur elements $(s_E)_{E \in \Irr(G_n)}$ for $\mathcal{H}(G_n)$ as calculated by Malle in \cite{Ma2}. We can also calculate the irreducible characters $(\chi_E)_{E \in \Irr(G_n)}$  of the split semisimple algebra $K(\textbf{{v}})\mathcal{H}(G_n)$ as the traces of its irreducible representations. Thus, using Formula \eqref{needed}, we can define the (conjectural) canonical symmetrising trace 
\begin{equation}\label{tildetau}
\widetilde{\tau} := \sum_{E \in \Irr(G_n)}\frac{1}{s_E}\chi_E
\end{equation}
on  $K(\textbf{{v}})\mathcal{H}(G_n)$.  By showing that $\widetilde{\tau}(b)=\delta_{1b}$ for all $b \in \mathcal{B}_n$, we establish that $\tau=\widetilde{\tau}$ on $\mathcal{H}(G_n)$.  Hence,
 $A=(\widetilde{\tau}(bb'))_{b,b' \in \mathcal{B}_n}$.
 We already know that $A$ is symmetric, because $\widetilde{\tau}$ is a linear combination of characters, and thus a trace function. Therefore, we only need to calculate the entries in and over the diagonal.
 
 We have decided to include the GAP3 algorithm in the file for each group that can be found on the project's webpage \cite{Chaw}, in case some reader wants to try to run it on their own.

\subsection{The extra condition}\label{extra con section}

Our algorithms so far have allowed us to determine the entries of the matrix   $A:=(\tau(bb'))_{b,b' \in \mathcal{B}_n}$, where $\tau(b):=\delta_{1b}$ for $b \in \mathcal{B}_n$, for all $n=4,\ldots,8$. In each of these cases, we can easily check that the matrix $A$ is symmetric and that its determinant is a unit in $\Z[\textbf{{u}},\textbf{{u}}^{-1}]$ {(we will give the exact values of the determinants in Section \ref{Results})}. Thus, $\tau$ is a symmetrising trace on $\mathcal{H}(G_n)$. Further, given that 
there is a section $G_n \rightarrow \boldsymbol{G}_n$, $w \mapsto \boldsymbol{w}$ of $G_n$ in $B(G_n)$ such that  $\mathcal{B}_n = \{T_{\boldsymbol{w}}\,|\,\boldsymbol{w} \in \boldsymbol{G}_n\}$, $\tau$ specialises to the canonical symmetrising trace on the group algebra of $G_n$ when $u_{\mathcal{C},j} \mapsto \zeta_{e_\mathcal{C}}^j$. 
Therefore, in order for the BMM symmetrising trace conjecture to hold, it remains to  check Condition \eqref{extra con}. Since our symmetrising trace satisfies the lifting conjecture by Malle and Michel, it suffices to check Condition \eqref{extra con2}, or the equivalent Condition \eqref{extra con3}.

At first, we  decided to modify our GAP3 algorithm in order to verify Condition \eqref{extra con3} for the map $\widetilde{\tau}$ defined in \eqref{tildetau}.  This was before we realised that our program could not even establish that $\tau=\widetilde{\tau}$ for the groups $G_5$ and $G_7$.  However, in the cases of the  three  remaining groups, where we have checked that $\tau=\widetilde{\tau}$, it was easy to verify that Condition \eqref{extra con3} holds. This second GAP3 algorithm is also included in the file for each group that can be found on the project's webpage \cite{Chaw}. At this point, we would like to thank Jean Michel for sending us the function $m$ which computes the values of the integers $m_{\mathcal{C},j}^E$ appearing in the formula for $\omega_E(T_{\boldsymbol{\pi}})$ given by \cite[Equation (1.22)]{BMM}.

For the groups $G_5$ and $G_7$, we had to check by hand  Condition \eqref{extra con2} directly.  
It is extremely interesting that one of our inputs for the SAGE algorithm are the coefficients of the element $T_{\boldsymbol{\pi}}=z^{|Z(G_n)|}$ which appears in the extra condition for $\tau$. 

Given that $\mathcal{B}_n = \{T_{\boldsymbol{w}}\,|\,\boldsymbol{w} \in \boldsymbol{G}_n\}$, checking
Condition \eqref{extra con2} amounts to showing that
\begin{equation}\label{extra con4}
\tau\left(b^{-1}z^{|Z(G_n)|}\right)=0, \quad \text{ for all } b \in \mathcal{B}_n \setminus \{1\}.
\end{equation}
The fact that the elements of $\mathcal{B}_n$ are (mostly) of the form $z^kb_m$, for $b_m$ in the much smaller set $\mathcal{E}_n$ and $k \in \{0,1,\ldots,|Z(G_n)|-1\}$, made the calculations easier. 
We also had the right to use that $\tau$ is a trace function.
The calculations can be found in the appendix of this paper, and they confirm the validity of Condition \eqref{extra con2}. 

In the end, we decided to check by hand  Condition \eqref{extra con2} also for the groups $G_4$, $G_6$ and $G_8$. Since we already know that the condition is satisfied for these groups thanks to the (second) GAP3 algorithm, these calculations are omitted from this paper, but can be found on the project's webpage \cite{Chaw}.

\section{Results}\label{Results}

We summarise here the results of the programs based on the algorithms presented in the previous section, concluding with the main result of this paper, the proof of the BMM symmetrising trace conjecture for  groups $G_4$, $G_5$, $G_6$, $G_7$, $G_8$. We also discuss some further consequences regarding our bases; notably we obtain a verification of the BMR freeness conjecture for these groups.

\subsection{On the BMM symmetrising trace conjecture}
Let $n \in \{4,\ldots,8\}$ and let $\mathcal{B}_n$ be the basis for $\mathcal{H}(G_n)$ given in 
Section \ref{Algorithm}. Let $\tau : \mathcal{H}(G_n) \rightarrow \Z[\textbf{{u}},\textbf{{u}}^{-1}]$ be the linear map defined by $\tau(b):=\delta_{1b}$ for all $b \in \mathcal{B}_n$.
Set  $A:=(\tau(bb'))_{b,b' \in \mathcal{B}_n}$.

\begin{thm}
The matrix $A$ is symmetric and we have
$$
{\rm det}(A) = \left\{\begin{array}{ll}\smallbreak \smallbreak
 -c^{58}\, & \text{for $n=4$};\\ \smallbreak \smallbreak
-c^{288}f^{276}\,& \text{for $n=5$}; \\ \smallbreak \smallbreak
b^{252}e^{200}\, & \text{for $n=6$}; \\ \smallbreak \smallbreak
b^{792}e^{648}h^{582}\, & \text{for $n=7$};\\\smallbreak \smallbreak
d^{600}\, & \text{for $n=8$}.
\end{array}\right.
$$
In particular, the determinant of $A$ is a unit in  $\Z[\textbf{{\em u}},\textbf{{\em u}}^{-1}]$. Thus, $\tau$ is a symmetrising trace on $\mathcal{H}(G_n)$.
\end{thm}

Now, by the construction of  $\mathcal{B}_n$,
there exists a section $G_n \rightarrow \boldsymbol{G}_n$, $w \mapsto \boldsymbol{w}$ of $G_n$ in $B(G_n)$ such that  $\mathcal{B}_n = \{T_{\boldsymbol{w}}\,|\,\boldsymbol{w} \in \boldsymbol{G}_n\}$.
Hence, by its definition, $\tau$ specialises to the canonical symmetrising trace on the group algebra of $G_n$ when  $u_{\mathcal{C},j}\mapsto \zeta_{e_\mathcal{C}}^{j}$. Moreover, following our discussion in Subsection \ref{extra con section}, $\tau$ satisfies Condition  \eqref{extra con2}, that is,
$$\tau(T_{\boldsymbol{w}^{-1}\boldsymbol{\pi}})=0, \quad \text{ for all } \boldsymbol{w} \in \boldsymbol{G}_n \setminus \{1\}.
$$
We obtain thus the main result of our paper:

\begin{thm}
Let $n \in \{4,\ldots,8\}$. The map $\tau$ is the canonical symmetrising trace on $\mathcal{H}(G_n)$.
In particular, the BMM symmetrising trace conjecture holds for $G_n$.
\end{thm}

In the process, we have also proved the following:
\begin{thm}
Let $n \in \{4,\ldots,8\}$. The lifting conjecture holds for $G_n$.
\end{thm}

\begin{rem}\rm
In Section \ref{Algorithm}, we saw that for $n \in \{5,6,8\}$, we took $\mathcal{B}_n$ to be the basis for $\mathcal{H}(G_n)$ given by the second author in \cite{Ch17}. We also explained in detail why we had to change the basis in the case of $G_7$. However, for $G_4$ we just chose $\mathcal{B}_4$ because all the generators appeared with positive exponents. Now, in the basis for $\mathcal{H}(G_4)$ given in \cite{Ch17}, the element $t^{-1}$ appears. Following the inverse Hecke relations in this case, we observe that $\tau(t^{-1}) =-bc^{-1}$. Given that the canonical symmetrising trace is unique, this implies that the change of basis was necessary for our approach.
\end{rem}

\subsection{On the BMR freeness conjecture}
Until now we have used the fact that the set $(\mathcal{B}_n)_{n=4,\ldots,8}$ is a basis for the generic Hecke algebra   $(\mathcal{H}(G_n))_{n=4,\ldots,8}$. We will now see that our \texttt{C++} program provides a proof of this fact.

Let $n \in \{4,5,6,8\}$.
We recall that our \texttt{C++} program has expressed $sb_j$ and $tb_j$ as $\Z[\textbf{{u}},\textbf{{u}}^{-1}]$-linear combinations of the elements of $\mathcal{B}_n$, for all $b_j \in \mathcal{B}_n$. This in fact allows us to express any product of the generators $s$ and $t$, and in turn any element of $\mathcal{H}(G_n)$, as a linear combination of the elements of $\mathcal{B}_n$.

\begin{exmp}\label{exmp G5}
\rm Let $n=5$. Using the notation of \S\ref{G5},
we have 
$$st^2s=st \sum_{l}\lambda^t_{2,l}b_l = 
s  \sum_{l}\lambda^t_{2,l} \sum_r  \lambda^t_{l,r}b_r = \sum_{l}  \lambda^t_{2,l} \sum_r \lambda^t_{l,r} \sum_p \lambda^s_{r,p}b_p.$$
\end{exmp}

Thus, we can also deduce the following:

\begin{prop}
Let $n \in \{4,5,6,8\}$.
The set $\mathcal{B}_n$ is a spanning set for $\mathcal{H}(G_n)$ as a $\Z[\textbf{{\em u}},\textbf{{\em u}}^{-1}]$-module.
\end{prop}

 Now, let $n=7$. The difference here is that our \texttt{C++} program has expressed $tb_j$ and $ub_j$ as linear combinations of the  elements of $\mathcal{B}_7$, for all $b_j \in \mathcal{B}_7$, but not $sb_j$. If we had these linear combinations, we could also express any element
 of $\mathcal{H}(G_7)$ as a linear combination of the elements of $\mathcal{B}_7$.
 We will show here that we do have them (following the notation of \S \ref{G7}):
 \begin{itemize}
 	\item If {$j\in\{1,\dots,132\}\setminus\{10,11,12\}$}, then $sb_j=(stu)u^{-1}t^{-1}b_j=u^{-1}t^{-1}zb_j=u^{-1}t^{-1}b_{12+j}$.
 	We now apply consecutively the inverse Hecke relations and we write $sb_j$ as a linear combination of elements of the form $u^{\alpha}t^{\beta}b_{12+j}$, with $\alpha,\beta \in\{0,1,2\}$. We continue now as in Example \ref{exmp G5} for $G_5$. \smallbreak
 		\item If {$j\in\{10,11,12\}$}, then $b_j=tu^2t^r$ for $r\in\{0,1,2\}$. We have $sb_j=(stu)ut^r=zut^r$, which is equal to the element $b_{14}$ for $r=0$, $b_{18}$ for $r=1$, and $b_{20}$ for $r=2$. \smallbreak
 	\item If {$j\in\{133,\dots, 144\}$}, then $sb_j=ab_j+bs^{-1}b_j$. It is thus enough to express $s^{-1}b_j$ as a linear combination of the elements of $\mathcal{B}_7$. We have 
	$s^{-1}b_j=(s^{-1}u^{-1}t^{-1})tub_j=tuz^{-1}b_j=tub_{j-12} =\sum_{l,p} \lambda^u_{j-12,l}  \lambda^t_{l,p}b_p$.\smallbreak
 \end{itemize}

Hence, we can also deduce the following:

\begin{prop}
The set $\mathcal{B}_7$ is a spanning set for $\mathcal{H}(G_7)$ as a $\Z[\textbf{{\em u}},\textbf{{\em u}}^{-1}]$-module.
\end{prop}

Since $|\mathcal{B}_n| = |G_n|$ for all $n \in \{4,\ldots,8\}$, Theorem \ref{mention} implies the following:

\begin{thm}\label{Bn basis}
Let $n \in \{4,\ldots,8\}$. The set $\mathcal{B}_n$ is a basis for $\mathcal{H}(G_n)$ as a $\Z[\textbf{{\em u}},\textbf{{\em u}}^{-1}]$-module. In particular, the BMR freeness conjecture holds for $G_n$.
\end{thm}

\begin{rem}\rm One could  argue that we have just obtained a computerised proof of the BMR freeness conjecture for groups $G_4, G_5, G_6, G_7, G_8$. This is indeed the case for $G_4$. However, the ``special cases'' incorporated in the \texttt{C++} algorithm for the other groups are the exact calculations made  by hand by the second author for the proof of the BMR freeness conjecture in \cite{Ch17}.
\end{rem}

\section{Appendix}\label{Appendix}
\subsection{The special cases for $G_6$}
In this section we discuss the special cases appearing in the \texttt{C++} algorithm (see \S \ref{C++}). These special cases follow the calculations in   \cite[Appendix]{Ch17} and can be found for each group on the project's webpage \cite{Chaw}.
We prove here the special cases for the group $G_6$.\\\\
\textbf{Case 1:} For every $k\in\{1,2,3\}$, $l_1,l_2\in\{0,1\}$, and $m_1,m_2,n_1,n_2\in\{0,1,2\}$ we have:\\
\phantom{Case 3:aa} $t^{m_1}s^{l_1}t^{n_1}z^kt^{m_2}s^{l_2}t^{n_2}=z^k t^{m_1}s^{l_1}t^{n_1+m_2}s^{l_2}t^{n_2}$.\\\\
\textbf{Case 2:} For every  $k,m,n\in\{0,1,2\}$ we have: 
$z^kt^{m}stst^{n}=z^{k+1}t^{m-1}s^{-1}t^{n-1}$.\\\\
\textbf{Case 3:} For every $k\in\{1,2,3\}$, and $m,n\in\{0,1,2\}$ we have:\\
\phantom{Case 3:} 
$\begin{array}[t]{lcl}
z^kt^{m}st^2st^{n}&=&cz^kt^mstst^n+
dz^kt^ms^2t^n+
aez^kt^{m-1}st^n+
abez^kt^ms^{-1}t^{n-1}+
eb^2z^{k-1}t^{m+1}st^{n+1}.
\end{array}$\\\\
\textbf{Case 4:} For every $m,n\in\{0,1,2\}$ we have:\\
\phantom{Case 4a:}$\begin{array}[t]{lcl}
z^3t^{m}stst^{n}&=&az^3t^{m+1}st^n+
abz^3t^ms^{-1}t^{n+1}+
cb^2z^3t^{m}s^{-2}t^n+
db^2z^2t^{m+1}st^{n+1}+\\&&
+ab^2ez^2t^{m+n+1}+
b^3ezt^{m+1}st^2st^{n+1}.
\end{array}$\\\\
\textbf{Case 5:} For every $m,n\in\{0,1,2\}$ we have:\\
\phantom{Case 4a:}
$\begin{array}[t]{lcl}
t^{m}st^2st^{n}&=&
-ab^{-2}zt^{m}stst^{n-2}+
b^{-2}e^{-1}z^2t^{m-1}s^{-1}tst^{n-1}-b^{-2}ce^{-1}z^2t^{m+n-2}
-
\\&&-ab^{-2}de^{-1}zt^mstst^{n-1}
-b^{-1}de^{-1}zt^mst^{n}-ab^{-1}zt^{m+n-1}.
\end{array}$\\\\
\textbf{Case 6:}
$\begin{array}[t]{lcl}
z^4&=&az^3tstst+bcz^3tst+bdz^2tsts^2tst+
abez^2ts^2tst+ab^2ez^2tst^2+b^3ecz^2s^{-1}tst+\\&& +b^3edzts^2t^2st+
ab^3e^2zst^2st+ab^4e^2zts^{-1}tst+b^5e^2t^2st^3st.
\end{array}$\\\\

The first case is straightforward, since $z$ is central. However, if we do not use the commutativity of the centre as a special case, the program takes a lot of time in order to commute $z$ with any element. In the case of $G_8$, this took more than one hour of calculations just for a single element!

Case 2 uses only the definition of $z=tststs$. Case 3 uses first the equivalent positive Hecke relation $t^2=ct+d+et^{-1}$:
$z^kt^mst^2st^n=cz^kt^mstst^n+
dz^kt^ms^2t^n+ez^kt^mst^{-1}st^n$. We apply now the equivalent positive  Hecke relation 
$s=a+bs^{-1}$  to the element $ez^kt^mst^{-1}st^n$ (twice) and we obtain: $ez^kt^mst^{-1}st^n=eaz^kt^{m-1}st^n+ebz^kt^ms^{-1}t^{-1}st^n=aez^kt^{m-1}st^n+
abez^kt^ms^{-1}t^{n-1}+eb^2z^kt^ms^{-1}t^{-1}s^{-1}t^n$. The result follows from the definition of $z$, which yields that the last term is equal to $eb^2 z^{k-1}t^{m+1}st^{n+1}$. 

We now prove Case 4. We use again the relation $s=a+bs^{-1}$, and we have:
$z^3t^{m}stst^{n}=az^3t^{m+1}st^n+bz^3t^{m}s^{-1}tst^{n}= az^3t^{m+1}st^n+abz^3t^{m}s^{-1}t^{n+1}+b^2z^3t^ms^{-1}ts^{-1}t^n$. We deal now with the last term of the sum. We use the relation $t=c+dt^{-1}+et^{-2}$ and we have:
$b^2z^3t^ms^{-1}ts^{-1}t^n=cb^2z^3t^ms^{-2}t^n+db^2z^3t^ms^{-1}t^{-1}s^{-1}t^n+eb^2z^3t^ms^{-1}t^{-2}s^{-1}t^n$. From the definition of $z$ we have that $db^2z^3t^ms^{-1}t^{-1}s^{-1}t^n=db^2z^2t^{m+1}st^{n+1}$ and $eb^2z^3t^ms^{-1}t^{-2}s^{-1}t^n=eb^2z^2t^{m+1}stst^{-1}s^{-1}t^n$. For the latter, we use again the relation $s=a+bs^{-1}$ and we obtain: $eb^2z^2t^{m+1}stst^{-1}s^{-1}t^n=ab^2ez^2t^{m+n+1}+b^3ez^2t^{m+1}sts^{-1}t^{-1}s^{-1}t^n$. The result follows again from the definition of $z$. We prove similarly Case 6.

Finally we prove Case 5. We underline the elements we encounter in Case 5 and we omit them from further calculations. Moreover, we write in bold the generator to which we apply the inverse Hecke relation. We have:
$$\begin{array}[t]{lcl}
	t^mst^2st^n&=&zt^mst\mathbf{s^{-1}}t^{-1}s^{-1}t^{n-1}\\
		&=&b^{-1}zt^mstst^{-1}\mathbf{s^{-1}}t^{n-1}-\underline{ab^{-1}zt^{m+n-1}}\\
			&=&b^{-2}zt^msts\mathbf{t^{-1}}st^{n-1}-\underline{ab^{-2}zt^mstst^{n-2}}\\
				&=&b^{-2}e^{-1}zt^mstst^2st^{n-1}-b^{-2}e^{-1}czt^mststst^{n-1}-b^{-2}e^{-1}dzt^msts^2t^{n-1}\\
				&=&\underline{b^{-2}e^{-1}z^2t^{m-1}s^{-1}tst^{n-1}}-\underline{b^{-2}ce^{-1}z^2t^{m+n-2}}-b^{-2}de^{-1}zt^mst(as+b)t^{n-1}\\
				&=&\underline{-ab^{-2}de^{-1}zt^mstst^{n-1}}-\underline{b^{-1}de^{-1}zt^mst^n}.
\end{array}$$

\subsection{The extra condition for $G_5$} We follow the notation of \S \ref{G5}. 
 Recall that $z=stst=tsts$ and $|Z(G_5)|=6$. In order to prove Condition \eqref{extra con4}, we will	 write $\tau(z^6b_i^{-1})$ for all ${b_i \in \mathcal{B}_5 \setminus \{1\}}$ as a $\mathbb{Z}[a,b,c^{\pm 1}, d,e,f^{\pm 1}]$-linear combination of elements of the form $\tau(b_jb_l)$ with ${b_j,b_l\in \mathcal{B}_5}$.
	Therefore, we can show that $\tau(z^6b_i^{-1})=0$ for all ${b_i\in \mathcal{B}_5 \setminus \{1\}}$ using the entries of the matrix $A$, which are computed by the \texttt{C++}/SAGE program.
	
	Let $b_i \in \mathcal{B}_5 \setminus \{1\}$. 
	We first consider the case where $i >12$.
	We write $b_i$ as $z^kb_m$, for $k\in\{1,\dots,5\}$ and $b_m\in \mathcal{E}_5$. Since $\tau$ is a trace function, we have that $\tau(z^{6-k}b_m^{-1})=\tau(z^{6-k}s^{p_1}t^{p_2})$, where $p_1\in\{-2,-1,0\}$ and $p_2\in\{-2,-1,0,1\}$. Using now the inverse Hecke relations, we can write 
	$\tau(z^{6-k}s^{p_1}t^{p_2})$ as a  $\mathbb{Z}[a,b,c^{\pm 1}, d,e,f^{\pm 1}]$-linear combination of elements of the form $\tau(z^{6-k}s^{q_1}t^{q_2})$, with $q_1,q_2\in\{0,1,2\}$. Since $k\in\{1,\dots,5\}$, we have that $z^{6-k}s^{q_1}t^{q_2}\in \mathcal{B}_5\setminus\{1\}$ and, hence, $\tau(z^{6-k}b_m^{-1})=0$. We now consider the case where $i \leq 12$. We have:\\
	\\{\small
	 $\begin{array}[t]{lcl}\tau(z^6b_2^{-1})&=&\tau(z^5tst)=\tau(z^5st^2)=\tau(b_{68})=0.\smallbreak\smallbreak\\
\tau(z^6b_3^{-1})&=&\tau(z^4tst^2st)=\tau(z^4st^2st^2)=\tau(b_{56}b_8)=0.\smallbreak\smallbreak\\
	\tau(z^6b_5^{-1})&=&\tau(z^4sts^2ts)=\tau(z^4s^2ts^2t)=\tau(b_{55}b_7)=0.\smallbreak\smallbreak\\
	\tau(z^6b_6^{-1})&=&\tau(z^5st)=\tau(b_{66})=0.\smallbreak\smallbreak\\
		\tau(z^6b_7^{-1})&=&\tau(z^5sts^{-1})=\tau(z^5t)=\tau(b_{64})=0.\smallbreak\smallbreak\\
		\tau(z^6b_8^{-1})&=&\tau(z^5tst^{-1})=\tau(z^5s)=\tau(b_{62})=0.\smallbreak\smallbreak\\
			\tau(z^6b_9^{-1})&=&\tau(z^5stst^{-1}s^{-2})=\tau(z^5tst^{-1}s^{-1})=\tau(z^4ts^2t)=\tau(z^4s^2t^2)=\tau(b_{57})=0.\smallbreak\smallbreak\\
			\tau(z^6b_{10}^{-1})&=&\tau(z^5tst^2)=\tau(z^5st^3)=\tau(b_{68}b_4)=0.\smallbreak\smallbreak\\
			\tau(z^6b_{11}^{-1})&=&\tau(z^5st^2)=\tau(b_{68})=0.\smallbreak\smallbreak\\
				\tau(z^6b_{12}^{-1})&=&\tau(z^5stst^{-1}s^{-1}t)=
				\tau(z^5tstst^{-1}s^{-1})=\tau(z^5st)=\tau(b_{66})=0.
	\end{array}$}\\ \\
	
	\subsection{The extra condition for $G_7$}
	We follow the notation of \S \ref{G7}.
	 Recall that $z=tus=ust=stu$ and $|Z(G_7)|=12$. In order to prove  Condition \eqref{extra con4}, we will	 write $\tau(z^{12}b_i^{-1})$ for all ${b_i \in \mathcal{B}_7 \setminus \{1\}}$ as a $\mathbb{Z}[a,b^{\pm 1}, c,d,e^{\pm 1},f,g,h^{\pm 1}]$-linear combination of elements of the form $\tau(b_jb_l)$ with ${b_j,b_l\in \mathcal{B}_7}$ and elements of the form $\tau(z^{12}b_r^{-1})$ with $r \neq i$ that have already been calculated.
	Therefore, we can show that $\tau(z^{12}b_i^{-1})=0$ for all ${b_i\in \mathcal{B}_7 \setminus \{1\}}$ using the entries of the matrix $A$, which are computed by the \texttt{C++}/SAGE program, along with some inductive arguments.

	Let $b_i\in \mathcal{B}_7 \setminus \{1\}$. We first consider the case where $i >12$.
	We write $b_i$ as $z^kb_m$, for $k\in\{1,\dots,11\}$ and $b_m\in \mathcal{E}_7$. Since $\tau$ is a trace function, we have that $\tau(z^{12-k}b_m^{-1})=\tau(z^{12-k}u^{p_1}t^{p_2})$, where $p_1\in\{-2,-1,0,1\}$ and $p_2\in\{-3,-2,-1,0\}$. Using now the inverse Hecke relations, we can write 
	$\tau(z^{12-k}u^{p_1}t^{p_2})$ as a $\mathbb{Z}[a,b^{\pm 1}, c,d,e^{\pm 1},f,g,h^{\pm 1}]$-linear combination of elements of the form $\tau(z^{12-k}u^{q_1}t^{q_2})$, with $q_1,q_2\in\{0,1,2\}$. Since $k\in\{1,\dots,11\}$, we have that $z^{12-k}u^{q_1}t^{q_2}\in \mathcal{B}_7\setminus\{1\}$ and, hence, $\tau(z^{12-k}b_m^{-1})=0$. We now consider the case where $i \leq 12$. We have:\\
	\\
	 ${\small
	 \begin{array}[t]{lcl}\tau(z^{12}b_2^{-1})&=&\tau(z^{11}st)=a\tau(z^{11}t)+b\tau(z^{11}s^{-1}t)
		=a\tau(b_{136})+b\tau(z^{10}tut)=b\tau(z^{10}ut^2)=b\tau(b_{128})=0.
	\smallbreak\smallbreak\\
	\tau(z^{12}b_3^{-1})&=&h^{-1}\tau(z^{12}u)-fh^{-1}\tau(z^{12})-gh^{-1}\tau(z^{12}u^{-1})=h^{-1}\tau(b_{13}b_{134})-
		fh^{-1}\tau(b_{73}b_{73})-gh^{-1}\tau(z^{12}b_2^{-1})=
		\smallbreak\smallbreak\\&=&
		h^{-1}(b^6e^4fh^4)-fh^{-1}(b^6e^4h^4)   =0.\smallbreak\smallbreak\\
\tau(z^{12}b_4^{-1})&=&\tau(z^{11}us)=a\tau(z^{11}u)+b\tau(z^{11}us^{-1})=a\tau(b_{134})+
b\tau(z^{10}utu)=b\tau(z^{10}u^2t)=b\tau(b_{127})=0.\smallbreak\smallbreak\\
\tau(z^{12}b_5^{-1})&=&e^{-1}\tau(z^{12}t)-ce^{-1}\tau(z^{12})-de^{-1}\tau(z^{12}t^{-1})=
e^{-1}\tau(b_{13}b_{136})-ce^{-1}\tau(b_{73}b_{73})-de^{-1}\tau(z^{12}b_4^{-1})=
\smallbreak\smallbreak\\&=&
e^{-1}(b^6ce^4h^4)-ce^{-1}(b^6e^4h^4)=0.\smallbreak\smallbreak\\
\tau(z^{12}b_6^{-1})&=&\tau(z^{11}s)=a\tau(z^{11})+b\tau(z^{11}s^{-1})=a\tau(b_{133})+b\tau(z^{10}tu)=b\tau(z^{10}ut)=b\tau(b_{126})=0.\smallbreak\smallbreak\\
\tau(z^{12}b_7^{-1})&=&\tau(z^{12}u^{-2}t^{-1})=\tau(z^{11}u^{-1}s)=\tau(z^{11}su^{-1})=a\tau(z^{11}u^{-1})+b\tau(z^{11}s^{-1}u^{-1})=\smallbreak\smallbreak\\&=&a\tau(z^{12}b_{14}^{-1})+b\tau(z^{10}t)= b\tau(b_{124})=0.\smallbreak\smallbreak\\
\tau(z^{12}b_8^{-1})&=&\tau(z^{12}u^{-1}t^{-2})=\tau(z^{11}st^{-1})=\tau(z^{11}t^{-1}s)=
a\tau(z^{11}t^{-1})+b\tau(z^{11}t^{-1}s^{-1})=\smallbreak\smallbreak\\&=&a\tau(z^{12}b_{16}^{-1}) +b\tau(z^{10}u)=b\tau(b_{122})=0.
\smallbreak\smallbreak\\
\tau(z^{12}b_9^{-1})&=&\tau(z^{11}u^{-1}st^{-1})=a\tau(z^{11}u^{-1}t^{-1})+b\tau(z^{11}u^{-1}s^{-1}t^{-1})=a\tau(z^{11}t^{-1}u^{-1})+b\tau(z^{10}u^{-1}tut^{-1})=\smallbreak\smallbreak\\&=&
a\tau(z^{12}b_{18}^{-1})+b\tau(z^{10}t^{-1}u^{-1}tu)=b\tau(z^9t^{-1}st^2u)=ab\tau(z^9tu)+b^2\tau(z^9t^{-1}s^{-1}t^2u)=\smallbreak\smallbreak\\&=&
ab\tau(z^9ut)+b^2\tau(z^8ut^2u)=ab\tau(b_{114})+b^2\tau(z^8u^2t^2)=b^2\tau(b_{105})=0.
\smallbreak\smallbreak\\

\tau(z^{12}b_{10}^{-1})&=&\tau(z^{12}u^{-2}t^{-1})=\tau(z^{12}t^{-1}u^{-2})=\tau(z^{12}b_7^{-1})=0.\smallbreak\smallbreak\\
\tau(z^{12}b_{11}^{-1})&=&\tau(z^{12}t^{-1}u^{-2}t^{-1})=\tau(z^{12}t^{-2}u^{-2})=\tau(b^{12}b_9^{-1})=0.\smallbreak\smallbreak\\
\tau(z^{12}b_{12}^{-1})&=&\tau(z^{12}t^{-2}u^{-2}t^{-1})=\tau(z^{12}t^{-3}u^{-2})=
e^{-1}\tau(z^{12}u^{-2})-ce^{-1}\tau(z^{12}t^{-1}u^{-2})-de^{-1}\tau(z^{12}t^{-2}u^{-2})=\smallbreak\smallbreak\\&=&e^{-1}\tau(z^{12}b_3^{-1})-ce^{-1}\tau(z^{12}b_{7}^{-1})-de^{-1}\tau(z^{12}b_{9}^{-1})=0.\smallbreak\smallbreak\\
	\end{array}}$

\end{document}